\documentclass[12pt]{article}
\usepackage{geometry}
\usepackage[latin1]{inputenc}
\usepackage{aeguill}
\usepackage[english]{babel}
\usepackage{amsfonts}
\usepackage{amssymb}
\usepackage{dsfont}
\usepackage{xspace}
\usepackage{amsmath}
\usepackage{amsthm}
\usepackage{indentfirst}
\usepackage[all]{xy}
\usepackage{hyperref}
\usepackage[OT2,T1]{fontenc}
\allowdisplaybreaks

\AtBeginDocument{
  \mathchardef\mathcomma\mathcode`\,
  \mathcode`\,="8000
} {\catcode`,=\active
  \gdef,{\mathcomma\discretionary{}{}{}}
}

\newtheorem{thm}{Theorem}[section]
\newtheorem{prop}[thm]{Proposition}

\newtheorem{lem}[thm]{Lemma}
\newtheorem{cor}[thm]{Corollary}

\newenvironment{rem}{\paragraph{Remark}}{\qed\vspace{1em}}
\newenvironment{dem}{\paragraph{Proof}}{\qed\vspace{1em}}

\newcommand{\N}{\mathrm{N}}
\newcommand{\Z}{\mathbb{Z}}
\newcommand{\Q}{\mathbb{Q}}
\newcommand{\R}{\mathbb{R}}
\newcommand{\C}{\mathbb{C}}
\newcommand{\Cp}{\C_{p}}

\newcommand{\Hom}{\mathrm{Hom}}
\newcommand{\Qp}{\Q_{p}}
\newcommand{\Zp}{\Z_{p}}

\newcommand{\A}{\mathbb{A}}
\newcommand{\Ax}{\A^{\times}}
\newcommand{\Af}{\A_{f}}
\newcommand{\Afx}{\Af^{\times}}

\newcommand{\AEx}{\A_{E}^{\times}}
\newcommand{\OEv}{\OF_{E,v}}
\newcommand{\OEE}{\OF_{E}}
\newcommand{\Ex}{E^{\times}}
\newcommand{\Fx}{F^{\times}}
\newcommand{\Fv}{F_{v}}
\newcommand{\OF}{\mathcal{O}}
\newcommand{\Ov}{\mathcal{O}_{v}}
\newcommand{\Cx}{\C^{\times}}
\newcommand{\GL}{\mathrm{GL}_{2}}
\newcommand{\PGL}{\mathrm{PGL}_{2}}
\newcommand{\SO}{\mathrm{SO}_{2}}
\newcommand{\Hn}{\mathcal{H}^{n}}
\newcommand{\dF}{d_{F}}

\newcommand{\T}{\mathrm{T}}
\newcommand{\TN}{\mathrm{T}_{N}}
\newcommand{\V}{\mathrm{V}}
\newcommand{\old}{\mathrm{old}}
\newcommand{\new}{\mathrm{new}}
\newcommand{\M}{\mathcal{M}}
\newcommand{\MS}{\mathcal{S}}
\newcommand{\ord}{\mathrm{ord}}
\newcommand{\val}{\mathrm{val}}
\newcommand{\PP}{\mathfrak{P}}
\newcommand{\DEF}{D_{E/F}}
\newcommand{\dEF}{d_{E/F}}
\newcommand{\GO}{\mathrm{GO}}
\newcommand{\Fvx}{\Fv^{\times}}
\newcommand{\sgn}{\mathrm{sgn}}
\newcommand{\Ev}{E_{v}}
\newcommand{\Evx}{\Ev^{\times}}
\newcommand{\AEE}{\A_{E}}
\newcommand{\Ovx}{\mathcal{O}_{v}^{\times}}
\newcommand{\OEx}{\mathcal{O}_{E}^{\times}}
\newcommand{\OFx}{\mathcal{O}^{\times}}
\newcommand{\dE}{d_{E}}
\newcommand{\AEf}{\A_{E,f}}
\newcommand{\AEfx}{\AEf^{\times}}
\newcommand{\Gal}{\mathrm{Gal}}
\newcommand{\mTheta}{\mathbf{\Theta}}
\newcommand{\mE}{\mathbf{E}}
\newcommand{\Bx}{B^{\times}}
\newcommand{\F}{\mathrm{F}}

\DeclareMathOperator{\Pic}{Pic}

\DeclareMathOperator{\Ext}{Ext}

\DeclareMathOperator{\Jac}{Jac}

\DeclareMathOperator{\rec}{rec}

\DeclareMathOperator{\Real}{Re}

\DeclareMathOperator{\Cl}{Cl}

\DeclareMathOperator{\Frob}{Frob}

\DeclareMathOperator{\End}{End}

\DeclareMathOperator{\tr}{tr}
\DeclareSymbolFont{cyrletters}{OT2}{wncyr}{m}{n}
\DeclareMathSymbol{\Sha}{\mathalpha}{cyrletters}{"58}

\title{$p$-adic Gross-Zagier Formula for Heegner Points on Shimura
Curves over Totally Real Fields}

\author{Li \textsc{Ma}}

\date{30 September 2014}

\begin{document}
\maketitle

  \section*{Abstract}
  
The main result of this text is a generalization of Perrin-Riou's p-adic Gross-Zagier formula to the case of Shimura curves over totally real fields. Let $F$ be a totally real field. Let $f$ be a Hilbert modular form over $F$ of parallel weight $2$, which is a new form and is ordinary at $p$. Let $E$ be a totally imaginary quadratic extension of $F$ of discriminant prime to $p$ and to the conductor of $f$. We may construct a $p$-adic $L$ function that interpolates special values of the complex $L$ functions associated to $f$, $E$ and finite order Hecke characters of $E$. The $p$-adic Gross-Zagier formula relates the central derivative of this $p$-adic $L$ function to the $p$-adic height of a Heegner divisor on a certain Shimura curve.

The strategy of the proof is close to that of the original work of Perrin-Riou. In the analytic part, we construct the analytic kernel via adelic computations; in the geometric part, we decompose the geometric kernel into two parts: places outside $p$ and places dividing $p$. For places outside $p$, the $p$-adic heights are essentially intersection numbers and are computed in works of S. Zhang, and it turns out that this part is closely related to the analytic kernel. For places dividing $p$, we use the method in the work of J. Nekov\'a\v r to show that the contribution of this part is zero.

\section*{Résumé}

Le résultat principal de ce texte est une généralisation de la formule de Gross-Zagier $p$-adique de Perrin-Riou au cas de courbes de Shimura sur les corps totalement réels. Soit $F$ un corps totalement réel. Soit $f$ une forme modulaire de Hilbert sur $F$ de poids parallel $2$, qui est une forme nouvelle et est ordinaire en $p$. Soit $E$ une extension quadratique totalement imaginaire de $F$ de discriminant premier à $p$ et au conducteur de $f$. On peut construire une fonction $L$ $p$-adique qui interpole valeurs spéciales de la fonction $L$ complexe associée à $f$, $E$ et caractères de Hecke d'ordre fini de $E$. La formule $p$-adique de Gross-Zagier relie la dérivée centrale de cette fonction $L$ $p$-adique à la hauteur d'un divisor de Heegner sur une certaine courbe de Shimura.

La stratégie de la preuve est proche de celle du travail original de Perrin-Riou. Dans la partie analytique, on construit le noyau analytique par calculs adéliques; dans la partie géométrique, on décompose le noyau géométrique en deux parties: places hors de $p$ et places divisant $p$. Pour les places hors de $p$, les hauteurs $p$-adiques sont essentiellement des nombres d'intersection et sont calculées dans les travaux de S. Zhang, et il s'avère que cette partie est bien liée au noyau analytique. Pour les places divisant $p$, on utilise la méthode dans le travail de J. Nekov\'a\v r pour montrer que la contribution de cette partie est nulle.

\tableofcontents

\clearpage

\section{Introduction}

The original Gross-Zagier formula, proved in \cite{GrossZagier},
expresses the central derivative of the $L$ function associated to a
new form of weight $2$ over $\Q$ and an imaginary quadratic field as
the (real-valued) height of a Heegner divisor on the modular
curve.\\

Later in \cite{Perrin-Riou87}, Perrin-Riou attached a $p$-adic $L$
function to a $p$-ordinary new form of weight $2$ and an imaginary
quadratic field, which interpolates special values of the complex
$L$ functions. She then proved that the central derivative of this
$p$-adic $L$ function can be expressed as the $p$-adic height of a
Heegner divisor. This was generalized by Nekov\'a\v{r}
\cite{Nekovar95} to the case of higher weight modular forms over
$\Q$.\\

On the other hand, the original (complex, or real) Gross-Zagier
formula has been generalized by Zhang \cite{Zhang01} (and subsequent
works) to Hilbert modular forms over any totally real field. In this
paper we would like to prove a $p$-adic version of Zhang's formula,
which is a generalization of Perrin-Riou's work to the case of any
totally real field.\\

Let $p$ be a fixed prime number. Fix injections
$\overline{\Q}\hookrightarrow\overline{\Qp}$ and
$\overline{\Q}\hookrightarrow\C$. Let $F$ be a totally real number
field, $\OF=\OF_{F}$ its integer ring and $N$ an ideal of $\OF$
prime to $p$. Let $f$ be a new form of level $N$, parallel weight
$2$ and trivial central character. The form $f$ admits a
$q$-expansion:
$$
f=\sum\limits_{0\neq I\subseteq\OF}a(I,f)q^I.
$$
We assume that the form $f$ is ordinary at $p$, i.e. the coefficient
$a(p\OF,f)$ has $p$-adic absolute value $1$.\\

Let $E$ be a totally imaginary quadratic extension of $F$ with
relative discriminant $D=D_{E/F}$, such that all prime ideals
dividing $2Np$ split in $E$. Then for any finite order character
$\chi:\AEx/\Ex\rightarrow\Cx$, the complex $L$ function is given by:
$$
L(s,f,\chi)=L(2s-1,\epsilon\chi_F)D(s,f,\chi),$$
where
$\epsilon:\A_{F}^{\times}/\Fx\rightarrow\Cx$ is the quadratic
character associated to the extension $E/F$, $\chi_F$ is the
restriction of $\chi$ to $F$, and $D(s,f,\chi)$ is:
$$
\sum\limits_{I\subseteq{\OEE}}a(\N_{E/F}(I),f)\chi(I)
\left|I\right|_\infty^{-s}
$$
(We use $\left|I\right|_\infty$ to denote the absolute norm of
$I$).\\

We are only interested in characters $\chi$ that factor through the
quotient $G=\AEx/\Ex\prod\limits_{v\nmid p}\OEx$, which is the
Galois group $\Gal(E^{(p)}/E)$ with $E^{(p)}$ the maximal unramified
outside $p$ abelian extension of $E$. Our first result is that there
exists a $p$-adic pseudo-measure $\mathbf{K}_f$ on the group $G$,
with possible poles only on characters $\chi$ such that $\chi_F$ is
the $p$-cyclotomic character, and satisfies the interpolation
property: for any finite order character $\chi:G\rightarrow\Cx$, the
integral $\int\limits_{G}\chi\mathbf{K}_f$ is, up to some simple
factor, equal to $L(1,f,\chi^{-1})$.\\

Let $\mathfrak{X}$ be the rigid analytic space
$\Hom_{\mathrm{cont}}(G,\Cp^{\times})$, then the pseudo-measure
$\mathbf{K}_f$ induces a meromorphic function $L_p$ on the space
$\mathfrak{X}$ via: $L_p(x)=\int\limits_{G}x\mathbf{K}_f$, which is
analytic on the anti-cyclotomic line
$\{x\in\mathfrak{X}:x=\overline{x}\}$. Thus we may take the
derivative of the function $L_p$ on the trivial character and in the
direction of the cyclotomic character $\xi_E$,
$L'_{p,\xi_E}(1)=\frac{d}{ds}L_p(\xi_E^s)\big|_{s=0}$.\\

We would like to relate this derivative to the $p$-adic height of a
Heegner divisor on a certain Shimura curve. We make the assumption
that $\epsilon(N)=(-1)^{[F:\Q]-1}$. Fix an archimedean place $\tau$
of $F$. The assumption $\epsilon(N)=(-1)^{[F:\Q]-1}$ ensures that
there is a quaternion algebra $B$ over $F$, which ramifies exactly
at all infinite places different from $\tau$ and all finite places
$v$ such that $\epsilon_v(N)=-1$. Since every prime ideal dividing
$N$ splits in the field $E$, we may embed $E$ into $B$ and thus view
$E$ as a sub-algebra of $B$.\\

Let $R$ be an order of $B$ of type $(N,E)$, i.e. an order that
contains $\OEE$ and has discriminant $N$. The existence of such an
order is proved in \cite{Zhang01} section 1.5. We then have a
Shimura curve $X$, whose complex points are given by:
$$
X(\C)=B^{\times}\backslash\mathcal{H}^{\pm}\times
\widehat{B}^{\times}/\widehat{F}^{\times}\widehat{R}^{\times}.
$$
Here we let $B^{\times}$ act on the Poincar\'e double half plane
$\mathcal{H}^{\pm}$ via the fixed place $\tau$ and a fixed
isomorphism $B_\tau\simeq M_2(\R)$.\\

Unlike the classical case, there is no "cusp" for the curve $X$ if
the field $F$ is not equal to $\Q$, which we always assume. But
Zhang \cite{Zhang01} defined a canonical divisor class
$\xi\in\Pic(X)\otimes\Q$, the "Hodge class", which has degree $1$ on
every connected component of $X$. This replaces the cusp $\infty$ in
the classical case, and allows us to define a map
$\phi:X\rightarrow\Jac(X)\otimes\Q$ by sending a point $y$ to the
class of $y-\xi$.\\

Under the complex description, there are \textbf{CM points} on the
curve $X$, which are represented by pairs
$(z_0,b)\in\mathcal{H}^{\pm}\times\widehat{B}^{\times}$, where $z_0$
is the unique fixed point of $\Ex$ in $\mathcal{H}^{\pm}$. Shimura's
theory shows that these points are algebraic and defined over the
maximal abelian extension of $E$.\\

A CM point that is defined over the Hilbert class field $H$ of $E$
is called a "Heegner point". The group $\Gal(H/E)$ acts simply
transitively on the set of Heegner points. Let $z$ be the divisor
class
$\frac{1}{\#(\OEx/\OFx)}\sum\limits_{x}\phi(x)\in\Jac(X)\otimes\Q$,
where the sum runs over all Heegner points. This is the trace (or
the average) of any Heegner point. Let $z_f$ be the $f$-isotypic
part of $z$ via the Jacquet-Langlands correspondence.\\

The general theory of $p$-adic height pairings of Zarhin and
Nekov\'a\v{r} (\cite{Zarhin90}, \cite{Nekovar93}) constructs a
$p$-adic height pairing $\langle\cdot,\cdot\rangle$ on
$\Jac(X)(E)\times\Jac(X)(E)$. Let
$\mathfrak{P}_1,\cdots\mathfrak{P}_l$ be all prime ideals of $\OF$
above $p$, and for each $i$, let $\alpha_i$ be the unique root of
the polynomial $X^2-a(\PP_i,f)X+\left|\PP_i\right|_{\infty}$ that is
a $p$-adic unit. Our main result can be stated as:

\begin{thm}
We have the following identity:
$$
L_{p,\xi_E}'(1)=\left(\prod\limits_{i}\frac{(\alpha_i-1)^3}
{(\alpha_i+1)(\alpha_i^2-\left|\PP_i\right|_{\infty})}\right)\langle
z_f,z_f\rangle.
$$
\end{thm}

\begin{rem}
In \cite{Disegni13}, D. Disegni obtained independently the same
result. It turns out that the analytic part of his work is almost
the same as ours (with slight difference in the construction of
Eisenstein series), but the geometric methods are different.
\end{rem}

Let us mention one application of this theorem. As in the classical
case, let $A$ be a modular elliptic curve over $F$ with conductor
$N$, parameterized by the curve $X$, i.e. we have a surjective
morphism $X\rightarrow A$, and hence a morphism
$\Jac(X)\rightarrow\Jac(A)=A$. We assume that $A$ has ordinary good
reduction at places above $p$, thus corresponds to a $p$-ordinary
new form $f$ of parallel weight $2$ and level $N$. Let $z_A$ be the
projection of the divisor $z$ on $A$, which is an $E$-rational point
on the curve $A$. We then have a $p$-adic $L$ function
$L_p(A_E,\cdot)$ associated to the curve $A$ and to the quadratic
extension $E$, which interpolating special values of the complex $L$
functions $L(A_E,\chi,1)$, and our formula becomes:
$$
L_{p,\xi_E}'(A_E,1)=C\langle z_A,z_A\rangle,
$$
where $C$ is a non-zero factor.

\begin{cor}
If the order of vanishing of the $p$-adic $L$ function
$L_p(A_E,\cdot)$ at the trivial character is exactly $1$, then the
group $\Z\cdot z_A$ has finite index in $A(E)$, and the
Tate-\v{S}afarevi\v{c} group $\Sha(A/E)$ is finite.
\end{cor}

\begin{dem}
This follows from a generalization of Kolyvagin's method in the
classical case, c.f. for example \cite{Nekovar07} Theorem 3.2.
\end{dem}

\section*{Acknowledgements}

I would like to express my deep gratitude to my thesis advisor, Jan
Nekov\'a\v r, for suggesting me doing this problem, for everything
he taught me and for his constant help during the work. This article
would never have been completed without him.

I would like to thank Jacques Tilouine and Jan Hendrik Bruinier for
carefully reading my thesis and for their valuable comments. I would
also like to thank Bernadette Perrin-Riou and Joseph Oesterl\'e for
agreeing to serve on my thesis committee.

I would like to thank Marie-France Vign\'eras for teaching me a lot
on automorphic representations during my m\'emoire, which become
useful in the current work.

I would like to thank Li Wenwei for his clarification on Weil
representations.

I would like to thank all my colleagues and friends in Paris for
their useful discussion and warm friendship. I would also like to
thank my teachers during my study at \'Ecole Normale Sup\'erieure
for their guidance and help.

I would like to thank my parents and my wife for supporting me all
the time.
\newpage

\paragraph{Notations and conventions}
We fix the following default settings throughout the text:\\
Fix injections $\overline{\Q}\hookrightarrow\C$ and
$\overline{\Q}\hookrightarrow\overline{\Qp}$;\\
$p$: an odd rational prime number;\\
$F$: a totally real field, with integer ring $\OF$;\\
$N$: an ideal of $F$ as in the introduction;\\
$\A=\A_{F}$: the ad\`ele ring of $F$;\\
$\dF$: the absolute different of $F$;\\
$\PP_{1},\cdots,\PP_{l}$: prime ideals of $\OF$ above $p$;\\
$\PP$: the product of all $\PP_{i}$;\\
$\val_{1},\cdots,\val_{l}$: corresponding valuations;\\
$E$: a totally imaginary quadratic extension of $F$ as in the
introduction, with relative discriminant $D=D_{E/F}$;\\
$\dEF$: the relative different of $E/F$;\\
$\xi:\Ax/\Fx\rightarrow\widehat{\Z}^{\times}$: the cyclotomic
character;\\
$\xi_p:\Ax/\Fx\rightarrow\Zp^{\times}$: the $p$-adic cyclotomic
character;\\
$\psi:\A/F\rightarrow\Cx$: the additive character defined by
$\psi=\psi_{\Q}\circ\tr_{F/\Q}$, with $\psi_{\Q}$ the usual additive
character on $\A_{\Q}/\Q$ such that $\psi_{\Q,\infty}(x)=e^{2\pi
ix}$;\\
If $I$ is a fractional ideal of $\OF$, denote by
$\left|I\right|_{\infty}$ the (absolute) norm of $I$, i.e. the
index of $I$ in $\OF$;\\
If $A$ is any $\Z$-module, write $\widehat{A}$ for
$A\otimes_{\Z}\widehat{\Z}$;\\
A subscript $v$ usually indicates
"local component at the place $v$", and a subscript $\infty$
indicates local component at all the infinite places, e.g.
$F_{\infty}$ stands for $F\otimes_{\Q}\R$.\\

If $N$ is an ideal of $\OF$, a function $\phi:\Ax/\Fx\rightarrow\C$
is called a \textbf{function mod $N$} if it factors through the
quotient space
$\Ax/\Fx(1+\widehat{N})^{\times}F_{\infty,>0}^{\times}$. A Hecke
character $\chi:\Ax/\Fx\rightarrow\Cx$ is called a \textbf{character
mod $N$} if it is a function mod $N$. If $\chi$ is a character mod
$N$, it is called \textbf{primitive} if it is not a character mod
$D$ for any proper divisor $D$ of $N$. In this case, the ideal $N$
is called the \textbf{conductor} of the character $\chi$.\\

If $\phi$ is a function mod $N$, we may define a function on the
integral ideals of $\OF$, denoted by $\phi_{[N]}$, such that:
$$
\phi_{[N]}(I)=
\begin{cases}
\phi(i), & \textrm{if }I\textrm{ is prime to }N\textrm{ and
}i\in\Afx\textrm{ satisfies }i\OF=I\textrm{ and }i_{v}=1\textrm{ for
all }v\mid N;\\
0, & \textrm{if }I\textrm{ is not prime to }N.
\end{cases}
$$

The Haar measure on $\A/F$ is normalized to be self-dual with
respect to the character $\psi$.\\

All Frobenius maps are by default the "geometric Frobenius", and the
reciprocity map of class field theory is normalized in the geometric
way.\\

\newpage

\centerline{\Large{Part I. Analytic part}}

\section{Hilbert Modular Forms}

\subsection{Definition and $q$-expansion}

Let $N$ be an ideal of $\OF$. Define the following open compact
subgroups of $\GL(\widehat{\OF})$:
\begin{eqnarray*}
K_{0}(N) &:=& \left\{
\begin{pmatrix}
a & b\\
c & d
\end{pmatrix}\in\GL(\widehat{\OF}):c\equiv0\mod N\right\},\\
K_{1}(N) &:=& \left\{
\begin{pmatrix}
a & b\\
c & d
\end{pmatrix}\in K_{0}(N):d\equiv1\mod N\right\}.
\end{eqnarray*}
Also define $K_{\infty}:=\SO(F_{\infty})$, the standard maximal
compact subgroup of $\GL(F_{\infty})$.\\

Let $\kappa$ be a positive integer. A \textbf{modular form of
(parallel) weight $\kappa$ for $K_1(N)$} is a measurable function
$\phi:\GL(\A)\rightarrow\C$ such that:
\begin{enumerate}
\item $\phi(\gamma g)=\phi(g)$ for any $\gamma\in\GL(F)$;
\item $\phi(gkr(\theta)a)=\phi(g)e^{i\kappa\theta}$ for any
$k\in K_{1}(N)$, any $\theta\in F_{\infty}$ and any $a$ in
$F_{\infty,>0}$, where for any $\theta\in F_{\infty}$, we write
$r(\theta)$ for the element $
\begin{pmatrix}
\cos\theta & \sin\theta\\
-\sin\theta & \cos\theta
\end{pmatrix}$ in $K_{\infty}$, and we write $e^{i\kappa\theta}$ for
the product $\prod\limits_{v\mid\infty}e^{i\kappa\theta_v}$;\\
\item $\phi$ is slowly increasing, i.e. for every real number $c>0$
and any compact subset $\Omega$ of $\GL(\A)$, there exist constants
$C$ and $M$ such that $\phi\left(
\begin{pmatrix}
a & \\
 & 1
\end{pmatrix}g\right)\leq C\left|a\right|^{M}$ for any $g\in\Omega$
and any $a\in\Ax$ with $\left|a\right|\geq c$.
\end{enumerate}

Let $\chi:\Ax\rightarrow\Cx$ be a character mod $N$, such that
$\chi_{v}(-1)=(-1)^{\kappa}$ for every infinite place $v$. A
\textbf{modular form of weight $\kappa$ for $K_0(N)$ with central
character $\chi$} is a modular form $\phi$ of weight $\kappa$ for
$K_1(N)$ such that $\phi(ag)=\chi(a)\phi(g)$ for any $a\in Z(\Ax)$
(scalar matrices).\\

If $\phi$ is a modular form, then it has a Fourier expansion:
$$
\phi(g)=W_{0}(g)+\sum\limits_{\alpha\in\Fx}W\left(
\begin{pmatrix}
\alpha & \\
 & 1
\end{pmatrix}g\right),
$$
where:
\begin{eqnarray*}
W_{0}(g) &=& \int\limits_{\A/F}\phi\left(
\begin{pmatrix}
1 & x\\
 & 1
\end{pmatrix}g\right)dx,\\
W(g) &=& \int\limits_{\A/F}\phi\left(
\begin{pmatrix}
1 & x\\
 & 1
\end{pmatrix}g\right)\psi(-x)dx.
\end{eqnarray*}
The functions $W_{0}$ and $W$ satisfy the same condition 2 as the
modular form $\phi$, and have the following property:
\begin{eqnarray*}
W_{0}\left(
\begin{pmatrix}
1 & x\\
 & 1
\end{pmatrix}g\right) &=& W_{0}(g),\\
W\left(
\begin{pmatrix}
1 & x\\
 & 1
\end{pmatrix}g\right) &=& W(g)\psi(x).
\end{eqnarray*}

If we have $W_{0}(g)=0$ for almost all $g$, then the modular form
$\phi$ is called \textbf{cuspidal}; if for any fixed elements
$x_{f}\in\A_{f}$ and $y_{f}\in\Afx$, the function
\begin{eqnarray*}
\Hn &\rightarrow& \C\\
x_{\infty}+iy_{\infty} &\mapsto& y_{\infty}^{-\kappa/2}\phi\left(
\begin{pmatrix}
y_{f} & x_{f}\\
 & 1
\end{pmatrix}
\begin{pmatrix}
y_{\infty} & x_{\infty}\\
 & 1
\end{pmatrix}\right)
\end{eqnarray*}
is holomorphic on $\Hn$, then the modular form $\phi$ is called
\textbf{holomorphic}. Here $y_{\infty}^{-\kappa/2}$ stands for
$\prod\limits_{v\mid\infty}y_{v}^{-\kappa/2}$.

\begin{prop}
Let $\phi$ be a holomorphic modular form for $K_{1}(N)$. Then there
exists a unique function on the set of non-zero ideals of $\OF$,
$I\mapsto a(I)$, and a function on the narrow ideal class group
$\Ax/\Fx F_{\infty,>0}^{\times}\widehat{\OF}^{\times}$, $y\mapsto
a_{0}(y)$, such that for any $x\in\A$ and $y\in\Ax$ with
$y_{\infty}>0$, we have:
$$
\phi
\begin{pmatrix}
y & x\\
 & 1
\end{pmatrix}=\left|y\right|^{\kappa/2}(a_{0}(y_{f}\dF)+\sum\limits_{\alpha>0}
a(\alpha y_{f}\dF)\psi(i\alpha y_{\infty})\psi(\alpha x)),
$$
where the sum ranges through totally positive elements $\alpha$ in
$F$, and $a(I)$ is understood to be zero if $I$ is not an integral
ideal. Furthermore, the modular form $\phi$ is determined by these
functions.
\end{prop}

If $\phi$ is a holomorphic form, we may formally write
$$
\phi=a_{0}+\sum\limits_{0\neq I\subseteq\OF}a(I)q^{I}.
$$
This is called the \textbf{$q$-expansion} of $\phi$. The form $\phi$
is cuspidal if and only if the function $a_{0}$ is zero. We
sometimes write $a(I,\phi)$ to emphasize the dependence on $\phi$.

\begin{dem}
The proof is the same as \cite{Zhang01} Proposition 3.1.2, but note
that the constant term depends on the narrow ideal class, which
seems to be an error of \textit{loc. cit}.
\end{dem}

The space of holomorphic modular forms (resp. cusp forms) of weight
$\kappa$, level $N$ and central character $\chi$ will be denoted
$\M_{\kappa}(N,\chi)$ (resp. $\MS_{\kappa}(N,\chi)$). When the
central character $\chi$ is trivial, we denote also by
$\M_{\kappa}(N)$ and $\MS_{\kappa}(N)$ the corresponding spaces.

\subsection{Hecke operators}\label{heckeoperator}

Let $M$ be a non-zero ideal of $\OF$. The set
$$
U(M):=\left\{
\begin{pmatrix}
a & b\\
c & d
\end{pmatrix}\in M_{2}(\widehat{\OF}):c\equiv0\mod
N,(ad-bc)\OF=M\right\}
$$
is left and right invariant by $K_0(N)$, and thus defines a Hecke
operator $\T(M)$ on the space $\M_{\kappa}(N,\chi)$:
$$
(\T(M)\phi)(g):=\N(M)^{\kappa/2-1}\int_{U(M)}\phi(gh)dh,
$$
where the Haar measure on $\GL(\Af)$ is normalized such that the set
$K_{0}(N)$ has volume $1$. Equivalently, if we write $U(M)$ as a
disjoint union $\bigsqcup\limits_{j}h_{j}K_{0}(N)$, then we have:
$$
(\T(M)\phi)(g):=\N(M)^{\kappa/2-1}\sum\limits_{j}\phi(gh_{j}).
$$

From the definition it is clear that the operators $\T_{M}$ are
multiplicative, i.e. if $M$ and $M'$ are coprime, then we have
$\T(MM')=\T(M)\T(M')$.

\begin{prop}
Let $\phi$ be a form in the space $\mathcal{M}_{\kappa}(N,\chi)$ and
$M$ be a non-zero ideal of $\OF$. If the $q$-expansions of the forms
$\phi$ and $\T(M)\phi$ are $a_{0}+\sum a(L)q^{L}$ and $b_{0}+\sum
b(L)q^{L}$, respectively, then we have:
\begin{eqnarray*}
b(L) &=&
\sum\limits_{J\mid(L,M)}\chi_{[N]}(J)\N(J)^{\kappa-1}a(LM/J^{2}),\\
b_{0}(L) &=& \sum\limits_{J\mid
M}\chi_{[N]}(J)\N(J)^{\kappa-1}a_{0}(LM/J^{2}).
\end{eqnarray*}
\end{prop}

\begin{dem}
This is the classical calculation of $q$-expansion of Hecke
operators. A proof may be found in \cite{Zhang01} Proposition 3.1.4.
\end{dem}

\begin{cor}
Let $P$ be a prime ideal of $\OF$. Then on the space
$\mathcal{M}_{\kappa}(N,\chi)$ we have, for any integer $m\geq2$:
$$
\T(P^{m})=\T(P^{m-1})\T(P)-\chi_{[N]}(P)\N(P)^{\kappa-1}\T(P^{m-2}),
$$
Thus all the Hecke operators $\T_{M}$ commute with each other, and
satisfy the following formal identity:
$$
\sum\limits_{M}\T(M)\N(M)^{-s}=\prod\limits_{P}(1-\T(P)\N(P)^{-s}+
\chi_{[N]}(P)\N(P)^{-2s+\kappa-1})^{-1}
$$
\end{cor}

\begin{rem}
This can also be deduced directly from the definition of the Hecke
operators (interpreted as operators on the Bruhat-Tits tree).
\end{rem}

We would also like to define the $\V$ operator. Let $N$ and $D$ be
ideals of $\OF$. Let $d\in\Afx$ be a generator of $D$. We then
define the operator $\V_{D}:\M_{\kappa}(N,\chi)\rightarrow
\M_{\kappa}(ND,\chi)$ by:
$\V_{D}(\phi)(g):=\N(D)^{-\kappa/2}\phi\left(g
\begin{pmatrix}
d^{-1} & \\
 & 1
\end{pmatrix}\right)$. An easy calculation shows that if the form
$\phi$ has $q$-expansion $a_{0}(I)+\sum\limits_{0\neq
I\subseteq\OF}a(I)q^{I}$, then the form $\V_{D}(\phi)$ has
$q$-expansion $a_{0}(D^{-1}I)+\sum\limits_{0\neq
I\subseteq\OF}a(I)q^{DI}$.

\subsection{Peterson inner product}

Let $\phi_{1}$ and $\phi_{2}$ be two modular forms of same weight
and same central character. Then we may define their Peterson
product as:
$$
(\phi_{1},\phi_{2}):=\int \limits_{\GL(F)\backslash\GL(\A)/Z(\A)}
\overline{\phi_{1}(g)}\phi_{2}(g)dg,
$$
where $dg$ is the right $\GL(\A)$-invariant Borel measure on the
space $\GL(F)\backslash\GL(\A)/Z(\A)$ (cf. \cite{Platonov94} Theorem
3.17 for the existence and uniqueness up to a positive scalar). The
measure is induced by the measure on $\GL(\A)/Z(\A)$ normalized as
follows:
$$
\int\limits_{\GL(\A)/Z(\A)}f(g)dg:=
\int\limits_{\Ax}\int\limits_{\A}\int\limits_{K_{0}(1)}f\left(
\begin{pmatrix}
y & x\\
 & 1
\end{pmatrix}k\right)\left|y\right|^{-1}dkdxd^{\times}y.
$$

Now fix a level $N$, a weight $\kappa$ and a character $\chi$. The
Peterson inner product is then a positive definite Hermitian form on
the space $\mathcal{S}_{\kappa}(K_{0}(N),\chi)$, making it a Hilbert
space.

\subsection{Fricke and Atkin-Lehner involutions and adjoints
of Hecke operators}

Still fix $N$, $\kappa$ and $\chi$. We define an application:
\begin{eqnarray*}
w_{N}:\M_{\kappa}(N,\chi) &\rightarrow& \M_{\kappa}(N,\chi^{-1})\\
\phi &\mapsto& \left(g\mapsto\phi\left(g
\begin{pmatrix}
 & -1\\
n_{f} &
\end{pmatrix}\right)\chi^{-1}(\det g)\right),
\end{eqnarray*}
where $n_{f}\in\Afx$ is any element such that $n_{f}\OF=N$. It is
easy to see that this operator is well-defined. It is called the
\textbf{Fricke involution}. The word "involution" comes from
the fact that $w_{N}^{2}$ is the identity map.\\

More generally, suppose that the ideal $N$ factors as a product
$N=MD$ with $(M,D)=1$, and suppose that the character $\chi$ also
factors as $\chi=\chi_{M}\chi_{D}$, with $\chi_{M}$ (resp.
$\chi_{D}$) a character mod $M$ (resp. mod $D$). Then we may
similarly define the (partial) \textbf{Atkin-Lehner involutions}
$w_{M}$ by: $(w_{M}\phi)(g):=\phi\left(g
\begin{pmatrix}
 & -1\\
m_{f} &
\end{pmatrix}\right)\chi_{M}^{-1}(\det g)$, with
$m_{f}\in\Afx$ a generator of $M$. This operator takes values in the
space $\M_{\kappa}(N,\chi_{M}^{-1}\chi_{D})$.\\

In particular, in the case that $\chi$ is trivial, we see that for
any factor $M$ of $N$ such that $(M,N/M)=1$, the operator
$w_{M}:\M_{\kappa}(N)\rightarrow \M_{\kappa}(N)$ is defined.\\

The Atkin-Lehner involutions can be used to describe adjoints of
Hecke operators. On the Hilbert space $\MS_{\kappa}(N)$, let
$\T(M)^{*}$ denote the adjoint of the operator $\T(M)$. We then
have:

\begin{prop}\label{adjoint}
\begin{enumerate}
\item $T(M)^{*}=w_{N}^{-1}T(M)w_{N}$;
\item If $M$ is prime to $N$, then we have: $T(M)^{*}=T(M)$.
\end{enumerate}
\end{prop}

\begin{dem}
Since Hecke operators are commutative and generated by $T_{P}$ for
$P$ prime ideals of $\OF$, it suffices to prove the results for
$\T(P)$. This can be done by direct computation, remarking that if
$K_{0}(N)
\begin{pmatrix}
p & \\
 & 1
\end{pmatrix}K_{0}(N)$ has decomposition
$\bigsqcup\beta_{j}K_{0}(N)$, then $K_{0}(N)
\begin{pmatrix}
p^{-1} & \\
 & 1
\end{pmatrix}K_{0}(N)$ has decomposition $\bigsqcup
p^{-1}u_{N}^{-1}\beta_{j}u_{N}K_{1}(N)$, where $u_{N}=
\begin{pmatrix}
 & -1\\
n_f &
\end{pmatrix}$ with $n_f$ a generator of $N$. Here $p\in\Afx$
is a fixed generator of $P$.
\end{dem}

\subsection{Old forms and new forms}

Let $N$ be an ideal of $\OF$. Let $M$ be a proper divisor of $N$ and
write $D$ for the quotient $N/M$. If $\phi$ is a modular form in
$\MS_{\kappa}(M)$, then the form $\V_{D}(\phi)$ is in
$\MS_{\kappa}(N)$.\\

The subspace of $\MS_{\kappa}(N)$ generated by all the forms
$\V_{D}(\phi)$ (for $M$ running over all the proper divisors of $N$
and $\phi$ running over all the forms in $\MS_{\kappa}(M)$) is
called the space of \textbf{old forms}, denoted by
$\MS_{\kappa}(N)^{\old}$. Any form in this subspace will be called
an old form.\\

The subspace of \textbf{new forms} is the orthogonal complement of
the space of old forms with respect to the Peterson inner product,
denoted by $\MS_{\kappa}(N)^{\new}$.

\begin{prop}
The subspaces $\mathcal{S}_{\kappa}(K_{0}(N))^{\old}$ and
$\mathcal{S}_{\kappa}(K_{0}(N))^{\new}$ are stable under Hecke
operators.
\end{prop}

\begin{dem}
It suffices to prove that the old forms are stable under the Hecke
operators and their adjoints. By Proposition \ref{adjoint}, it
suffices to prove that the old forms are stable under the operators
$\T(M)$ and $w_{N}$. These are all easy to verify.
\end{dem}

\section{$p$-adic modular forms and Hida's theory}

\subsection{$p$-adic modular forms and Hecke operators}

We will define $p$-adic modular forms \textit{à la} Serre.\\

Let $R$ be any ring. The ring of formal $q$-expansions with
coefficients in $R$, denoted by $R[[q]]$, is the set of all formal
series of the form $a_{0}+\sum\limits_{0\neq
I\subseteq\OF}a(I)q^{I}$, where $a_{0}$ is a function on the narrow
ideal class group of $F$ and $a$ is a function on the set of
non-zero ideals of $\OF$, both taking values in $R$. The addition of
two series in $R[[q]]$ is the pointwise addition, and multiplication
is defined as if they were $q$-expansions of modular forms: if
$f=a_{0}+\sum\limits a(I)q^{I}$ and $g=b_{0}+\sum\limits b(I)q^{I}$
are two series in $R[[q]]$, their product $fg$ will be the formal
series $c_{0}+\sum\limits c(I)q^{I}$ such that:
\begin{eqnarray*}
c_{0}(I) &=& a_{0}(I)b_{0}(I),\\
c(I) &=& \sum\limits_{\substack{\alpha,\beta\in
I^{-1}\\\alpha\geq0,\beta\geq0\\\alpha+\beta=1}}a(\alpha I)b(\beta
I),
\end{eqnarray*}
where in the last equation, the value $a(\alpha I)$ for $\alpha=0$
is understood to be $a_{0}(I)$.\\

If $M$ is an ideal of $\OF$, $f$ is a holomorphic form for
$K_{0}(M)$ and $R$ is a sub-ring of $\C$ that contains all
coefficients of $f$, then we may identify $f$ with its
$q$-expansion, considered as an element of the ring
$R[[q]]$.\\

If $R=L$ is a fixed $p$-adic field, the modular forms lie in the
smaller ring of \textbf{bounded formal $q$-expansions}, $L\langle
q\rangle$, which is by definition the sub-ring of $L[[q]]$
consisting of those series whose coefficients $a$ and $a_{0}$ are
bounded functions. The ring $L\langle q\rangle$ is equipped with the
$\sup$ norm on its coefficients, which makes it an $L$-Banach
algebra.\\

Let $\kappa$ be a positive integer. For $M$ a non-zero ideal of
$\OF$, let $\M_{\kappa}(M,\Q)$ denote the intersection of
$\M_{\kappa}(M)$ with $\Q[[q]]$. The complex vector space
$\M_{\kappa}(M)$ is finite dimensional and has a basis consisting of
forms having rational $q$-expansions, i.e. we have:
$\M_{\kappa}(M,\Q)\otimes\C=
\M_{\kappa}(M)$.\\

Define $\M_{\kappa}(M,L)$ to be the image of
$\M_{\kappa}(M,\Q)\otimes L$ in the space $L[[q]]$.
Thus we know that this image lies in the subspace $L\langle q\rangle$.\\

Denote by $\M_{\kappa}(Mp^{\infty},L)$ the closure (with respect to
the topology of $L\langle q\rangle$ induced by its norm) of the
union of all the spaces $\M_{\kappa}(Mp^{r},L)$ for all $r>0$. This
is then an $L$-Banach space.\\

If $\chi$ is any finite order character of conductor dividing
$Mp^{\infty}$, we define in the same way the space
$\M_{\kappa}(Mp^{\infty},\chi,L)$, just by enlarging $\Q$ to a
number field that contains all values of $\chi$.\\

Finally, we define $\MS_{\kappa}(Mp^{\infty},L)$ or
$\MS_{\kappa}(Mp^{\infty},\chi,L)$ to be those $p$-adic modular
forms with zero constant terms.\\

It turns out that Hecke operators on the space of (complex) modular
forms induce operators on the space of $p$-adic modular forms.\\

In fact, for all $r>0$, the ideals $Mp^{r}$ always have the same
prime divisors, so the Hecke operators always have the same effects
on the coefficients of the $q$ expansion, and the effects are
obviously continuous with respect to the $\sup$ norm. Thus the
$p$-adic Hecke operators can be defined via these actions on
$q$-expansions.

\subsection{Hida's theory}

Let $N$ be an ideal of $\OF$ prime to $p$.\\

A $p$-adic modular form in $\M_{\kappa}(Np^{\infty},L)$ is called
"ordinary" if the coefficients $a(\mathfrak{P}_{i})$ are $p$-adic
units for each $\mathfrak{P}_{i}$ dividing $p$. Hida's theory
(\cite{Hida85}) shows that the ordinary forms in
$\M_{\kappa}(Np^{\infty},L)$ are essentially in
$\M_{\kappa}(Np,L)$.\\

More precisely, the operator
$e_{\ord}=\lim\limits_{k\rightarrow\infty}\T(p)^{k!}$ acts on the
space $\M_{\kappa}(Np^{\infty},L)$, and takes values in
$\M_{\kappa}(N\mathfrak{P},L)$, where $\mathfrak{P}$ is the ideal
$\prod\limits_{\mathfrak{P}_{i}\mid p}\mathfrak{P}_{i}$. If $f$ is a
Hecke eigenform, then we have:
$$
e_{\ord}(f)=
\begin{cases}
f, & \textrm{ if $f$ is ordinary};\\
0, & \textrm{ else}.
\end{cases}
$$

Now let $f$ be the modular form in the introduction, i.e. a cuspidal
new form of weight $2$, of level $N$ prime to $p$ and of trivial
central character, which is a Hecke eigenform and is ordinary at
$p$. When we regard $f$ as a form in $\M_{2}(Np^{\infty},L)$,
however, it is no longer an eigenvector for the Hecke operators
dividing $p$, because they are different from the corresponding
level $N$ operators. Thus we define its $p$-stabilization as
follows. Let $\alpha_{i}$ (resp. $\beta_{i}$) be the unique root of
the polynomial
$X^{2}-a(\mathfrak{P}_{i},f)X+\left|\mathfrak{P_{i}}\right|_\infty$
that is (resp. is not) a $p$-adic unit. Then the form $f_{0}$
defined by:
$$
f_{0}=\sum\limits_{J\subseteq\{1,\cdots,l\}}(-1)^{\#J}
\left(\prod_{j\in J}\beta_{j}\right)\V_{\prod\limits_{j\in
J}\mathfrak{P}_{j}}(f)
$$
is the unique ordinary form in $\MS_{2}(Np^{\infty}, L)$ which is an
eigenvector for all Hecke operators, and has the same eigenvalues as
$f$ for Hecke operators with index prime to $p$. It actually is a
form in
$\MS_{2}(N_{0})$, where $N_0$ is the ideal $N\mathfrak{P}$.\\

By results of Hida, the line $L\cdot f_{0}$ is a direct summand of
the ordinary part, regarded as modules over the Hecke algebra
generated by all Hecke operators. Let $e_{f_{0}}$ be the
corresponding projector. We then define a linear form $l_{f}$ on
$\M_{2}(Np^{\infty},L)$ by:
$l_{f}(g):=a(1,e_{f_{0}}e_{\mathrm{ord}}(g))$.\\

This linear form then "interpolates" the Peterson inner product with
$f_{0}$ in the following sense. If $h$ is a holomorphic modular
form, let $h^{\star}$ be the form such that
$h^{\star}(g):=\overline{h\left(g
\begin{pmatrix}
-1 & \\
 & 1
\end{pmatrix}\right)}$, which is another holomorphic modular form,
with Fourier coefficients $a(I,h^{\star})=\overline{a(I,h)}$ for any
ideal $I$. We then have:

\begin{prop}\label{hidaformula}
The linear form $l_{f}$ on $\mathcal{M}_{2}(K_{0}(Np^{\infty}),L)$
has the following interpolation property: if $P$ is a divisor of
$p^{\infty}$, and $g$ is a form in
$\mathcal{M}_{2}(K_{0}(N_{0}P),L)$ with algebraic Fourier
coefficients, then we have:
$$
l_{f}(g)=\left|P\right|_{\infty}\left(\prod\alpha_{i}^{-\val_{i}
(P)}\right)\frac{(w_{N_{0}P}(f_{0}^{\star}),g)}{(w_{N_{0}}
(f_{0}^{\star}),f_{0})}
$$
\end{prop}

\begin{dem}
This is similar to Hida's original proof. If
$\gamma=(\gamma_1,\cdots)$ is any array in
$\prod\left\{\alpha_{i},\beta_{i}\right\}$, let $f_{\gamma}$ be the
form
$\left[\prod\limits_{i}(1-\gamma_{i}V_{\mathfrak{P}_{i}})\right]f$.
Thus if $\gamma$ is the array such that $\gamma_i=\beta_i$ for all
$i$, then $f_{\gamma}$ is just the form $f_{0}$. Without confusion,
we denote this special array by $0$.\\

In the case $P=1$, the only thing to verify is that for any array
$\gamma\neq0$,  the product
$S_{\gamma}:=(w_{N_{0}}(f_{0}^{\star}),f_{\gamma})$ is zero. Suppose
without loss of generality that $\gamma_{1}=\alpha_{1}$, then this
is shown by the following calculations:
\begin{eqnarray*}
&& (w_{N_{0}}(f_{0}^{\star}),\T(\PP_{1})(f_{\gamma}))=
(w_{N_{0}}(f_{0}^{\star}),\beta_{1}(f_{\gamma})) =
\beta_{1}S_{\gamma};\\
&=& (\T(\PP_{1})^{*}w_{N_{0}}(f_{0}^{\star}),f_{\gamma})=
(w_{N_{0}}T(\PP_{1})(f_{0}^{\star}),f_{\gamma})=
(\overline{\alpha_{1}}w_{N_{0}}(f_{0}^{\star}),f_{\gamma})=
\alpha_{1}S_{\gamma}.
\end{eqnarray*}

For general $P$, the operator $\T(P)$ takes a form $g$ in
$\mathcal{M}_{2}(K_{0}(N_{0}P),L)$ to
$\mathcal{M}_{2}(K_{0}(N_{0}),L)$, and we may compute:
\begin{eqnarray*}
&& l_{f}(\T(P)g)=a(1,e_{f_{0}}e_{\mathrm{ord}}\T(P)g)=
a(1,\T(P)e_{f_{0}}e_{\mathrm{ord}}g)\\
&=& a\left(1,\left(\prod\alpha_{i}^{\val_{i}(P)}\right)
e_{f_{0}}e_{\mathrm{ord}}g\right)=
\left(\prod\alpha_{i}^{\val_{i}(P)}\right)l_{f}(g).
\end{eqnarray*}
Applying the case $P=1$ to the form $\T(P)g$, we get:
\begin{eqnarray*}
&& \left(\prod\alpha_{i}^{\val_{i}(P)}\right)l_{f}(g)=
l_{f}(\T(P)g)\\
&=& \frac{(w_{N_{0}}(f_{0}^{\star}),\T(P)(g))}
{(w_{N_{0}}(f_{0}^{\star}),f_{0})}=\frac{(\T_{N_{0}P}(P)^{*}w_{N_{0}}
(f_{0}^{\star}),g)}{(w_{N_{0}}(f_{0}^{\star}),f_{0})}.
\end{eqnarray*}
Here $\T_{N_{0}P}(P)^{*}$ is the dual operator on the space
$\M_{2}(N_{0}P)$, thus is equal to $w_{N_{0}P}\T(P)w_{N_{0}P}$.
Moreover, it is easy to see (from the corresponding matrix identity)
that $w_{N_{0}P}w_{N_{0}}=\left|P\right|_{\infty}V_{P}$, and that
$\T(P)V_{P}$ is identity, so the form
$\T(P)^{*}w_{N_{0}}(f_{0}^{\star})$ simplifies to:
$\left|P\right|_{\infty}w_{N_{0}P}f_{0}^{\star}$, and gives the
desired result.
\end{dem}

From the definition and the fact that the operator
$\T(\mathfrak{P}_i)$ commutes with $e_{f_0}$ and $e_{\mathrm{ord}}$,
it is clear that we have
$l_f(\T(\mathfrak{P}_i)(g))=\alpha_{i}l_f(g)$ for every $i$. As in
the classical case, we may compare "Peterson inner product with
$f_0$" with "Peterson inner product with $f$".

\begin{lem}\label{lemmaproduct}
If $g$ is a modular form in the space $\M_2(N)$ with algebraic
coefficients, then we have:
$$
\left(\prod\limits_i(1-\frac{\left|\mathfrak{P}_i\right|_\infty}
{\alpha_i^2})\right)l_f(g)=\frac{(f,g)}{(f,f)}.
$$
\end{lem}

\begin{dem}
This is true for all modular forms $g$ with complex coefficients, if
we replace the left hand side by the Peterson inner product
expression in Proposition \ref{hidaformula}. The proof is the same
as in the classical case of elliptic modular forms
\cite{Perrin-Riou87} Lemme 2.2 (the detailed calculations are in
\cite{Perrin-Riou88} lemme 27).
\end{dem}

\section{Complex Theta functions}

In this section, let $N\neq1$ be an ideal of $\OEE$. Let
$\chi:\AEx/\Ex(1+\widehat{N})^{\times}E_{\infty}^{\times}
\rightarrow\C$ be a character mod $N$. Using the method of
\cite{YuanZhangZhang}, we are going to attach to the character
$\chi$ a holomorphic modular form in the space
$M_{1}(\DEF\N_{E/F}(N),\chi)$, the theta function.

\subsection{Weil representations}

For simplicity, here we only recall the Weil representation over $F$
attached to the quadratic space $(E,\N)$. Also, instead of
introducing the action of the whole group $\GO(\A)$ (cf.
\cite{YuanZhangZhang}
section 2.1), we only introduce the action of $\AEx$.\\

Let $v$ be a place of $F$. We define Schwartz functions on the space
$E_{v}\times\Fvx$ as follows: if $v$ is a finite place, a Schwartz
function is a locally constant and compactly supported function; if
$v$ is infinite, a Schwartz function is of the form:
$$
(t,u)\mapsto(P_{1}(u\N(t))+\sgn(u)P_{2}(u\N(t)))
e^{-2\pi\left|u\right|\N(t)},
$$
with $P_{1}$, $P_{2}$ polynomials over $\Fv$, and $\sgn$ is the sign
of $u$.\\

For any place $v$, we let $\Evx$ act on the space $E_{v}\times\Fvx$
via: $e\cdot(t,u):=(e^{-1}t,\N(e)u)$.

\begin{prop}
There is a unique representation of $\GL(\Fv)\times\Evx$ on the
space of Schwartz functions on $E_{v}\times\Fvx$, such that for any
Schwartz function $\phi$, we have:
\begin{itemize}
\item $(e\phi)(t,u)=\phi(e^{-1}(t,u))$, for any $e\in\Evx$;
\item $\left(
\begin{pmatrix}
y & \\
 & y^{-1}
\end{pmatrix}\phi\right)(t,u)=\phi(yt,u)\left|y\right|_{v}\epsilon_{v}(y)$,
for any $y\in\Fvx$;
\item $\left(
\begin{pmatrix}
1 & \\
 & d
\end{pmatrix}\phi\right)(t,u)=\phi(t,d^{-1}u)
\left|d\right|_{v}^{-1/2}$, for any $d\in\Fvx$;
\item $\left(
\begin{pmatrix}
1 & x\\
 & 1
\end{pmatrix}\phi\right)(t,u)=\phi(t,u)\psi(xu\N(t))$, for any
$x\in\Fv$;
\item $\left(
\begin{pmatrix}
 & 1\\
-1 &
\end{pmatrix}\phi\right)(t,u)=\gamma_{v}\epsilon_{v}(u)
\left|u\right|_{v}\int\limits_{E_{v}}\phi(s,u)\psi(u\N(t,s))ds$,
\end{itemize}
where $\gamma_{v}$ is a fourth root of unity (the Weil index)
depending only on $E$ and $F$, and $\N(t,s)$ stands for
$\N(t+s)-\N(t)-\N(s)$.\\

We also have: $\prod\limits_{v}\gamma_{v}=1$ (cf. \cite{Weil64}
Proposition 5).
\end{prop}

\begin{rem}
From the above formulae, it is easy to deduce the following:
\begin{itemize}
\item $\left(
\begin{pmatrix}
y & \\
 & 1
\end{pmatrix}\phi\right)(t,u)=\phi(yt,y^{-1}u)\left|y\right|_{v}^{1/2}
\epsilon_{v}(y)$, for any $y\in\Fvx$;
\item $\left(
\begin{pmatrix}
z & \\
 & z
\end{pmatrix}\phi\right)(t,u)=(z\phi)(t,u)\epsilon_{v}(z)$,
for any $z\in\Fvx$;
\end{itemize}
\end{rem}

In the global case, define Schwartz functions on the space
$\AEE\times\Ax$ as finite linear combinations of products
$\prod\limits_{v}\phi_{v}$ of Schwartz functions $\phi_{v}$ on each
place $v$, where for almost all $v$, the function $\phi_{v}$ is
spherical, i.e. is the characteristic function of $\OEv\times\Ovx$.
Then the local Weil representations induce a global representation
of the group $\GL(\A)\times\AEx$ on the space of Schwartz functions
on $\AEE\times\Ax$.

\subsection{Theta functions}\label{Thetasection}

Let $\phi$ be a Schwartz function on $\AEE\times\Ax$. Then there is
a subgroup $U$ of finite index of $\OEx$, such that $\phi$ is
invariant under the action of $U$ (cf. \cite{YuanZhangZhang} section
3.1). Since the action of $\GL(\A)$ commutes with the action of
$\AEx$, we see that $g\phi$ is invariant under $U$ for any
$g\in\GL(\A)$.

\begin{prop}
Let $\phi$ be a Schwartz function on $\AEE\times\Ax$ and $U$ be a
subgroup of finite index of $\OEx$ such that $\phi$ is invariant
under $U$. For any $g\in\GL(A)$, define:
$$
\theta_{\phi}(g):=\sum\limits_{(t,u)\in
U\backslash(E\times\Fx)}(g\phi)(t,u),
$$
then the sum converges uniformly for $g$ in any compact subset of
$\GL(\A)$, and defines an automorphic form on $\GL(\A)$.
\end{prop}

\begin{dem}
The only problem is the convergence, which is proved in
\cite{YuanZhangZhang} section 3.1.
\end{dem}

\begin{rem}
If the group $U$ in the proposition is replaced by some smaller
group $U'$, then the resulting function $\theta_{\phi}$ will be
$[U:U']$ times bigger. Thus the notation $\theta_{\phi}$ is
well-defined up to a constant factor.
\end{rem}

Now write $x\mapsto\overline{x}$ for the action of the non-trivial
element of $\Gal(E/F)$. Let $N$ be an ideal of $\OEE$. Define a
Schwartz function $\rho$ as follows:
\begin{eqnarray*}
\rho_{f}(t,u) &:=&
\mathds{1}_{tu\in(1+\widehat{\overline{N}})}\mathds{1}_{u\in\dF
\widehat{\OF}^{\times}}\epsilon(u);\\
\rho_{\infty}(t,u) &:=& \mathds{1}_{u>0}e^{-2\pi u\N(t)}.
\end{eqnarray*}
The function $\rho$ (and thus all the functions $\sigma\rho$ for
$\sigma\in\AEx$) is invariant under action of the group
$U:=\OEx\bigcap(1+\widehat{N})^{\times}$, so we may define, for
every $\sigma\in\AEx$, the associated theta function:
$$
\theta_{\sigma}(g):=\sum\limits_{U\backslash(E\times\Fx)}
(g\sigma\rho)(t,u).
$$
This is the \textbf{(mod $N$) partial theta function } associated to
$\sigma\in\AEx$.

\begin{lem}
The function $\rho$ has the following properties:
\begin{itemize}
\item $\left(
\begin{pmatrix}
1 & x\\
 & 1
\end{pmatrix}\rho_{f}\right)(t,u)=\rho_{f}(t,u)$, for any
$x\in\widehat{\OF}$;
\item $\left(
\begin{pmatrix}
y & \\
 & 1
\end{pmatrix}\rho_{f}\right)(t,u)=
\rho_{f}(t,u)$, for any $y\in\widehat{\OF}^{\times}$;
\item $\left(
\begin{pmatrix}
1 & \\
a & 1
\end{pmatrix}\rho_{f}\right)(t,u)=\rho_{f}(t,u)$, for
any $a\in D_{E/F}\widehat{\N(N)}$;
\item $\left(
\begin{pmatrix}
\cos\alpha & \sin\alpha\\
-\sin\alpha & \cos\alpha
\end{pmatrix}\rho_{\infty}\right)(t,u)=
\rho_{\infty}(t,u)e^{i\alpha}$, for any $\alpha\in F_{\infty}$;
\item $(\tau\rho)(t,u)=\rho(t,u)$, for any
$\tau\in(1+\widehat{N})^{\times}E_{\infty}^{\times}$;
\end{itemize}
\end{lem}

\begin{dem}
Only the third assertion needs some explanation. Write the matrix
$\begin{pmatrix}
1 & \\
a & 1
\end{pmatrix}$ as $W^{-1}
\begin{pmatrix}
1 & -a\\
 & 1
\end{pmatrix}W$, where $W$ denotes the matrix $
\begin{pmatrix}
 & 1\\
-1 &
\end{pmatrix}$. The function
$W\rho_{f}$ can be evaluated as follows:
\begin{eqnarray*}
W\rho_{f}(t,u) &=&
\gamma_{f}\mathds{1}_{u\in\dF\widehat{\OF}^{\times}}
\left|\dF\right|_{f}^{1/2}\int
\limits_{u^{-1}(1+\widehat{\overline{N}})}\psi(u\N(t,s))ds\\
&=&
A(t,u)\int\limits_{\widehat{\overline{N}}}\psi_{E}(\overline{t}s)ds,
\end{eqnarray*}
where $A$ is a function. It follows that the above function is
supported on the compact set $\{(t,u):u\in\dF\Ovx,
t\in\widehat{\overline{(d_{E}N)^{-1}}}\}$ (notice that
$\overline{\dE}=\dE$). Since for any $a\in D_{E/F}\widehat{\N(N)}$,
$u\in\dF\Ovx$ and $t\in\widehat{\overline{(d_{E}N)^{-1}}}$ the
number $au\N(t)$ belongs to $\dF^{-1}\Ov$, we have $
\begin{pmatrix}
1 & -a\\
 & 1
\end{pmatrix}W\rho_{f}=W\rho_{f}$ for any $a\in
D_{E/F}\widehat{\N(N)}$. So finally we get:
$$
\begin{pmatrix}
1 & \\
a & 1
\end{pmatrix}\rho_{f}=W^{-1}W\rho_{f}=\rho_{f},
$$
for any $a\in D_{E/F}\widehat{\N(N)}$.
\end{dem}

\begin{cor}
For any $\sigma\in\AEfx$, the function $\theta_{\sigma}$ is a
holomorphic modular form of weight $1$ for $K_{1}(D_{E/F}\N(N))$.\\

We also have: $\theta_{\sigma e\tau}=\theta_{\sigma}$ for any
$e\in\Ex$ and any
$\tau\in(1+\widehat{N})^{\times}E_{\infty}^{\times}$, so we may
define $\theta_{\sigma}$ for
$\sigma\in\AEx/\Ex(1+\widehat{N})^{\times}E_{\infty}^{\times}$.
\end{cor}

\begin{dem}
The holomorphicity is easily checked from the form of
$\rho_{\infty}$. The other things are clear from the above lemma.
For example, to prove that $\theta_{\sigma}$ is right invariant
under $K_{1}(D_{E/F}\N(N))$, take any
$d\in(1+D_{E/F}\widehat{\N(N)})^{\times}$ and compute:
$$
\begin{pmatrix}
1 & \\
 & d
\end{pmatrix}\rho_{f}=
\begin{pmatrix}
d & \\
 & d
\end{pmatrix}\rho_{f}=(d\rho_{f})\epsilon(d)=\rho_{f},
$$
then remark that for any ideal $M$, the group $K_{1}(M)$ is
generated by elements of the form $
\begin{pmatrix}
1 & x\\
 & 1
\end{pmatrix}$, $
\begin{pmatrix}
y & \\
 & 1
\end{pmatrix}$, $
\begin{pmatrix}
1 & \\
 & d
\end{pmatrix}$, $
\begin{pmatrix}
1 & \\
a & 1
\end{pmatrix}$, with $x\in\widehat{\OF}$, $y\in\widehat{\OF}^{\times}$,
$d\in(1+\widehat{M})^{\times}$, $a\in\widehat{M}$.
\end{dem}

Now let $\chi:\AEx/\Ex(1+\widehat{N})^{\times}\rightarrow\Cx$ be a
character mod $N$. The \textbf{theta function} associated to the
character $\chi$ is defined by:
$$
\Theta(\chi):=\sum\limits_{\sigma\in\AEx/\Ex(1+\widehat{N})^{\times}
E_{\infty}^{\times}}\chi(\sigma^{-1})\theta_{\sigma}.
$$

\begin{prop}
The theta function associated to $\chi$ is a holomorphic modular
form in the space $\M_{1}(\DEF\N(N),\epsilon\chi_{F})$, where
$\chi_{F}$ denotes the restriction of $\chi$ to $\Ax$.
\end{prop}

\begin{dem}

 The only thing
left to verify is the central character. In fact, for any $z\in\Ax$,
we have:
$$
\theta_{\sigma}\left(
\begin{pmatrix}
z & \\
 & z
\end{pmatrix}g\right)=\sum\limits_{(t,u)}(g\sigma
z\rho)(t,u)\epsilon(z)=\theta_{\sigma z}(g)\epsilon(z),
$$
so that:
$$
\Theta\left(
\begin{pmatrix}
z & \\
 & z
\end{pmatrix}g\right)=\sum\limits_{\sigma}\chi(\sigma^{-1})
\theta_{\sigma z}(g)\epsilon(z)=\Theta(g)(\chi\epsilon)(z).
$$
\end{dem}

\subsection{$q$-expansion of the theta function}

Here we compute the $q$-expansions of the partial theta functions.\\

For any element
$\sigma\in\AEx/\Ex(1+\widehat{N})^{\times}E_{\infty}^{\times}$,
write $\mathds{1}_{\sigma}:\AEx/\Ex(1+\widehat{N})^{\times}
E_{\infty}^{\times}\rightarrow\C$ for the characteristic function of
$\sigma$.

\begin{prop}
Let $\sigma$ be an element of
$\AEx/\Ex(1+\widehat{N})^{\times}E_{\infty}^{\times}$. Let $a$ and
$a_{0}$ be the coefficients of the $q$-expansion of
$\theta_{\sigma}$. Let $I$ be any non-zero ideal of $\OEE$. Then we
have:
\begin{eqnarray*}
a(I) &=& \sum\limits_{\substack{J\subseteq\OEE\\\N(J)=I}}
\mathds{1}_{\sigma^{-1}}(J),\\
a_{0}(I) &=& 0.
\end{eqnarray*}
\end{prop}

\begin{dem}
Let $I$ be any non-zero ideal of $\OEE$, and take any $y\in\Afx$
such that $y\dF=I$. We have:
\begin{eqnarray*}
a(I) &=& \left|y\right|^{-1/2}\psi(-iy_{\infty})
\int\limits_{\A/F}\theta_{\sigma}
\begin{pmatrix}
y & x\\
 & 1
\end{pmatrix}\psi(-x)dx\\
&=& \left|y\right|^{-1/2}\psi(-iy_{\infty})\sum\limits_{U\backslash
(E\times F^{\times})}\left(
\begin{pmatrix}
y & \\
 & 1
\end{pmatrix}\sigma\rho\right)(t,u)\int\limits_{\A/F}
\psi((u\N(t)-1)x)dx\\
&=& \sum\limits_{\substack{U\backslash(E\times
F^{\times})\\u\N(t)=1}}\rho_{f}(y\sigma
t,y^{-1}\N(\sigma)^{-1}u)\epsilon(y)\\
&=&
\sum\limits_{t\in\Ex/U}\mathds{1}_{(t\sigma)^{-1}\in(1+\widehat{N})}
\mathds{1}_{\N(t\sigma)^{-1}\in y\dF\widehat{O}^{\times}}\\
&=& \sum\limits_{\substack{J\subseteq\OEE\\\N(J)=I}}
\mathds{1}_{\sigma^{-1}}(J).
\end{eqnarray*}

The same calculation for the constant term gives:
$$
a_{0}(I)=\sum\limits_{u\in\Fx/U}\mathds{1}_{u>0}
\mathds{1}_{0\in(1+\widehat{\overline{N}})}
\mathds{1}_{u\N(\sigma)^{-1}\in y\dF\widehat{O}^{\times}}.
$$
This vanishes under the hypothesis that $N$ is not equal to $\OEE$.
\end{dem}

\begin{cor}
Let
$\chi:\AEx/\Ex(1+\widehat{N})^{\times}E_{\infty}^{\times}\rightarrow\C$
be a character mod $N$. Then the $q$-expansion of the modular form
$\Theta_{\chi}$ is given by:
$$
\theta_{\phi}=\sum\limits_{0\neq
J\subseteq\OEE}\chi_{[N]}(J)q^{\N(J)}.
$$
\end{cor}

\section{Complex Eisenstein series}\label{sectionEisenstein}

In this section, let $\kappa$ be a positive integer, $N$ be an ideal
of $\OF$ and $\chi:\Ax/\Fx\rightarrow\Cx$ be a Hecke character mod
$N$ and having the same parity as $\kappa$, i.e. such that
$\chi_{v}(-1)=(-1)^{\kappa}$ for any infinite place $v$. We are
going to recall the construction of an Eisenstein series in the
space $\M_{\kappa}(N,\chi)$.

\subsection{The case of primitive characters}

We first consider the case when the character $\chi$ is of conductor
$N$. For any finite place $v$ of $F$, define a function $H_{v}$ on
$\GL(\Ov)$ as follows: if $v$ does not divide $N$, put $H_{v}=1$;
otherwise, for $k=
\begin{pmatrix}
u & v\\
w & t
\end{pmatrix}\in\GL(\Ov)$, put
$H_{v}(k)=\mathds{1}_{K_{0}(N)}(k)\chi_{v}(t)$. Also, for $v=\infty$
define $H_{\infty}$ as: $H_{\infty}(r(\theta))=e^{i\kappa\theta}$,
where $r(\theta)$ is the matrix $
\begin{pmatrix}
\cos\theta & \sin\theta\\
-\sin\theta & \cos\theta
\end{pmatrix}$. The product of $H_{v}$ for all places $v$ then
defines a function $H$ on the group $K_{0}(1)K_{\infty}$, i.e.
$H(g)=\prod\limits_{v}H_{v}(g_{v})$.\\

Let $s$ be a complex number such that $\Real(s)>1$. The function $H$
then extends to a function $H_{s}:\GL(\A)\rightarrow\C$ via the
Iwasawa decomposition: for any $g\in\GL(\A)$, write $g=
\begin{pmatrix}
a & b\\
 & d
\end{pmatrix}k$ with $k\in K_{0}(1)K_{\infty}$. We then define:
$$
H_{s}(g):=\left|\frac{a}{d}\right|^{s}\chi(d)H(k).
$$
It is easy to see that the function $H_{s}$ is well defined (because
of the hypothesis on the parity of $\chi$), and it also decomposes
as a product: $H_{s}(g)=\prod\limits_{v}H_{s,v}(g_{v})$.\\

Define the Eisenstein series $G_{s}:\GL(\A)\rightarrow\C$ and its
normalization $E_{s}$ as follows:
\begin{eqnarray*}
G_{s}(g) &:=& \sum\limits_{\gamma\in
B(F)\backslash\GL(F)}H_{s}(\gamma
g),\\
E_{s} &:=& C_{s}G_{s},
\end{eqnarray*}
where we write $C_{s}$ for the number $2^{-n}L_{N}(1-2s,\chi)$. The
two Eisenstein series are absolutely convergent if $\Real(s)>1$. An
easy calculation shows that in this case the function $G_{s}$ (and
hence the function $E_{s}$) is a modular form for $K_{0}(N)$ of
weight $\kappa$ and of character
$\chi$.\\

For any ideal $I$ of $\OF$ and any integer $m$, write
$\sigma_{\chi_{[N]},m}(I):=\sum\limits_{J\mid
I}\chi_{[N]}(J)\N(J)^{m}$.

\begin{prop}
Let $\chi$ be as above.
\begin{itemize}
\item The function $E_{s}$, originally defined only for those $s$
with $\Real(s)>1$, can be meromorphically continued (with respect to
$s$) to all $s\in\C$, and is holomorphic at the point $s=\kappa/2$.
\item When $(F,\kappa,\chi)\neq(\Q,2,1)$, the form $E_{\kappa/2}$ is
a holomorphic modular form. The coefficients of its $q$-expansion
are given by:
$$
a(I)=\sigma_{\chi_{[N]},\kappa-1}(I),\textrm{ for any }I\neq0.
$$
The constant term is given by:
$$
a_{0}(I)=2^{-n}L_{N}(1-\kappa,\chi),\textrm{ for any }I\neq0,
$$
except for the case that $\kappa=1$ and $\chi$ is unramified on all
finite places. In that case, the constant term is:
$$
a_{0}(I)=2^{-n}(L_{1}(0,\chi)+
(-1)^{n}\chi_{[1]}(I)L_{1}(0,\chi^{-1})).
$$
\end{itemize}
\end{prop}

\begin{dem}
The proof is the same as \cite{Zhang01} Proposition 3.5.2. Here we
briefly recall the main steps. Let $W_{s}$ and $W_{0,s}$ be the
Whittaker functions of the form $E_{s}$. Then for any $x\in\A$ and
any $y\in\Ax$ with $y_{\infty}>0$ we have:
$$
E_{s}
\begin{pmatrix}
y & x\\
 & 1
\end{pmatrix}=W_{0,s}
\begin{pmatrix}
y & \\
 & 1
\end{pmatrix}+\sum\limits_{\alpha\in\Fx}W_{s}
\begin{pmatrix}
\alpha y & \\
 & 1
\end{pmatrix}\psi(\alpha x).
$$
Note that a modular form is determined by its values on these
matrices. Thus it suffices to prove analytic continuation of the
Whittaker functions.\\

From the Bruhat decomposition:
$$
\GL(F)=B(F)\bigsqcup\left(\bigsqcup\limits_{u\in F}B(F)
\begin{pmatrix}
 & 1\\
1 & u
\end{pmatrix}\right)
$$
we get:
\begin{eqnarray*}
W_{s}
\begin{pmatrix}
y & \\
 & 1
\end{pmatrix} &=&
C_{s}\int\limits_{\A/F}\left(H_{s}
\begin{pmatrix}
y & x\\
 & 1
\end{pmatrix}+\sum\limits_{u\in F}H_{s}
\begin{pmatrix}
 & 1\\
y & u+x
\end{pmatrix}\right)\psi(-x)dx\\
&=& C_{s}\int\limits_{\A}H_{s}
\begin{pmatrix}
 & 1\\
y & x
\end{pmatrix}\psi(-x)dx\\
&=& C_{s}\left|y\right|^{-s}\chi(y)\int\limits_{\A}H_{s}
\begin{pmatrix}
 & 1\\
1 & xy^{-1}
\end{pmatrix}\psi(-x)dx\\
&=& C_{s}\left|y\right|^{1-s}\chi(y)\prod\limits_{v}V_{s,v}(y_{v}),
\end{eqnarray*}
where for each place $v$, the functions $V_{s,v}:\Fv\rightarrow\C$
is defined as:
$$
V_{s,v}(y):=\int\limits_{\Fv}H_{s,v}
\begin{pmatrix}
 & 1\\
1 & x
\end{pmatrix}\psi_{v}(-yx)dx.
$$

The same calculation for the function $W_{0,s}$ yields:
$$
W_{0,s}
\begin{pmatrix}
y & \\
 & 1
\end{pmatrix}=C_{s}(\left|y\right|^{s}+\left|y\right|^{1-s}
\chi(y)\prod\limits_{v}V_{s,v}(0)).
$$

It then suffices to compute the functions $V_{s,v}$ for each place
$v$. We refer to \cite{Zhang01} for the details.
\end{dem}

We write $E_{\kappa/2}(\chi)$ for the holomorphic form
$E_{\kappa/2}$ defined above.

\subsection{The case of non-primitive characters}

Now assume that the character $\chi$ is of conductor $D$, a divisor
of $N$. We would also like to construct a modular form
$E=E_{\kappa/2}(\chi)$ in the space $\M_{\kappa}(N,\chi)$ with
Fourier coefficients $a(I)=\sigma_{\chi_{[N]},\kappa-1}(I)$ and
$a_{0}(I)=2^{-n}L_{N}(1-\kappa,\chi)$.\\

There are several possible constructions. Let $E_{D}$ be the modular
form constructed in the above subsection, i.e. the Fourier
coefficients $b(I)$ of $E_{D}$ is equal to
$\sigma_{\chi_{[D]},\kappa-1}(I)$. For any prime ideal $Q$ dividing
$N$, denote by $F(Q)$ the operator
$\T(Q)-\chi_{[D]}(Q)\N(Q)^{\kappa-1}$ on the space
$\M_{\kappa}(N,\chi)$. We then define the form $E$ as
$(\prod\limits_{Q\mid N})E_{D}$. It is then easy to verify that this
form $E$ has the desired Fourier coefficients.\\

Equivalently, we may define a function
$H_{s}=\prod\limits_{v}H_{s,v}$ on the group $K_{0}(1)K_{\infty}$ as
follows: for places $v$ dividing $N$ but not dividing $D$ and for
any $k=
\begin{pmatrix}
u & v\\
w & t
\end{pmatrix}\in\GL(\Ov)$, put:
$$
H_{s,v}(k)=
\begin{cases}
1-\left|\pi_{v}\right|, & \textrm{if }\left|w\right|_{v}=1;\\
1-\chi(\pi_{v})\left|\pi_{v}\right|^{1-2s}, & \textrm{if
}\left|w\right|_{v}<1,
\end{cases}
$$
and for other places $v$, define $H_{s,v}(k)$ the same way as the
above subsection. The function $H_{s}$ then extends to a function on
$\GL(\A)$ via the formula:
$$
H_{s}\left(
\begin{pmatrix}
a & b\\
 & d
\end{pmatrix}k\right)=\left|\frac{a}{d}\right|^{s}\chi(d)H_{s}(k).
$$
The Eisenstein series $G_{s}(g):=\sum\limits_{\gamma\in
B(F)\backslash\GL(F)}H_{s}(\gamma g)$ and its normalization
$E_{s}:=2^{-n}L_{D}(1-2s,\chi)G_{s}$ then admit holomorphic
extensions to all $s$, and the form $E_{\kappa/2}$ is holomorphic.
The same calculation as the subsection above shows that this form
has the desired $q$-expansion. We will use this construction in the
calculations below.\\

\begin{rem}
There is another construction of this Eisenstein series in
\cite{Disegni13}, which is $E:=w_{N}E_{D}$, where $E_{D}$ is the
form constructed in the above subsection. This defines the same
Eisenstein series as ours, because an easy calculation shows that
their Fourier coefficients are the same.
\end{rem}

\section{Measures on the Galois group}

From now on, let $E$ be a totally imaginary quadratic extension of
$F$, with relative discriminant $D=\DEF$ and relative difference
$\dEF$.

We would like to construct (pseudo-)measures on the Galois group
$G=\Gal(E^{(p)}/E)$, where $E^{(p)}$ is the maximal abelian
extension of $E$ unramified outside places dividing $p$. By class
field theory, we may rewrite the group $G$ as: $G=\AEx/\overline{\Ex
E_{\infty}^{\times}\prod\limits_{v\nmid p}\OEv^{\times}}$, here the
bar means taking closure.\\

The natural morphism $\prod\limits_{v\mid p}\OEv^{\times}\rightarrow
G$ induces an exact sequence:
$$
1\rightarrow\overline{\OEx}\rightarrow\prod\limits_{v\mid
p}\OEv^{\times}\rightarrow G\rightarrow\Cl(E)\rightarrow1.
$$
Thus the group $G$ is isomorphic to a direct product of a finite
group and some copies of $\Zp$.\\

\subsection{Theta measure}

To define the Theta measure, we first define a measure on the group
$G$ with values in $L\langle q\rangle$, the space of bounded formal
$q$-expansions, as follows: for any continuous function
$\phi:G\rightarrow L$, put
$$
\mTheta(\phi)=\sum\limits_{0\neq
J\subseteq\OEE}\phi_{[p]}(J)q^{\N(J)}.
$$
It is obvious that this is a measure on the group $G$.\\

By results of section \ref{Thetasection}, we have
$\mTheta(\chi)=\Theta(\chi)$ for any finite order character
$\chi:G\rightarrow\Cx$, so that $\mTheta(\chi)$ is actually in the
subspace $S_1(Dp^{\infty},\epsilon\chi_{F},L)$.

\subsection{Eisenstein pseudo-measure}

We first define a pseudo-measure on the group
$G_{F}:=\Ax/\overline{\Fx F_{\infty,>0}^{x}\prod\limits_{v\nmid
p}(1+D\Ov)^{\times}}$ with values in $L\langle q\rangle$.\\

For an element $s\in F_{\infty}^{\times}$, define its sign as
$\sgn(s)=\prod\limits_{v\mid\infty}\sgn(s_{v})$. For any function
$\phi:\Ax/\Fx F_{\infty,>0}^{\times}\rightarrow\C$, define:
$$
\phi^{\#}(x):=2^{-n}\sum\limits_{s\in
F_{\infty}^{\times}/F_{\infty,>0}^{\times}}\sgn(s)\phi(sx).
$$
By result of Deligne-Ribet, the map sending a locally constant
function $\phi:G_{F}\rightarrow L$ to the value
$L_{p\DEF}(0,\phi^{\#})$ extends to a pseudo-measure on $G_{F}$,
which we still denote by $L_{p\DEF}(0,\phi^{\#})$. More precisely,
this pseudo-measure $\mu$ is an element in the total quotient ring
of $\Lambda_{G_{F}}$ (the Iwasawa algebra of $G_{F}$), such that for
any element $C$ of the group $G_{F}$, the element
$\mu^{C}:=(1-\xi_{p}(C)^{-1}C)\mu$ lies in the sub-algebra of
($\Qp$-)measures over $G_{F}$.\\

We then define the Eisenstein pseudo-measure on $G_{F}$ by:
$$
\mathbf{E}_{F}(\phi):=2^{-n}L_{pD}(0,\phi^{\#})+
\sum\limits_{I}\sigma_{\phi^{\#}_{[pD]}}(I)q^{I},
$$
where for any function $\varphi$ on the set of integral ideals of
$\OF$, we define: $\sigma_{\varphi}(I):=\sum\limits_{J\mid
I}\varphi(J)$.\\

By results in section \ref{sectionEisenstein}, for any finite odd
character $\chi:G_{F}\rightarrow\Cx$, we have:
$\mE_{F}(\chi)=E_{1/2}(\chi)$, so that $\mE_{F}(\chi)$ is actually
in the subspace $\M_{1}(Dp^{\infty}),L)$. Also it is obvious that
$\mE_{F}(\chi)=0$ if the character $\chi$ is not odd.\\

Now define a pseudo-measure $\mE$ on the group $G$ as follows:
$$
\int_{G}\phi\mE:=\int_{G_{F}}\epsilon(g)\phi_{F}(g^{-1})\mE_{F}(g).
$$
For an element $C$ of $G_{F}$, we also write $\mathbf{E}^{C}$ for
the measure $\mathbf{E}(1-\xi_{p}(C)^{-1}C)$. If
$\phi=\chi:G\rightarrow\Cx$ is a finite order character, then (since
the character $\epsilon\chi_{F}$ is automatically odd)
$\mathbf{E}(\phi)$ is a modular form of central character
$\epsilon\chi_{F}^{-1}$.

\subsection{The analytic kernel}

Define the convoluted pseudo-measure
$\widetilde{\mathbf{K}}:=\mathbf{\Theta}(\V_{N}\mathbf{E})$. (And
$\widetilde{\mathbf{K}}^{C}$ means the measure
$\mathbf{\Theta}(\V_{N}\mathbf{E}^{C})$.)\\

If $\phi=\chi$ is a finite order character, then we have:
$$
\int\limits_{G}\chi\widetilde{\mathbf{K}}=\int\limits_{G\times
G}\chi(xy)\mathbf{\Theta}x(\V_{N}\mathbf{E})y=\mathbf{\Theta}(\chi)
\V_{N}\mathbf{E}(\chi),
$$
which is a modular form in the space $\mathcal{M}_{2}(NDp^{r})$ for
sufficiently large $r$. Since any locally constant function is a
linear combination of finite order characters, we see that the
pseudo-measure $\widetilde{\mathbf{K}}$ takes values in the $p$-adic
Banach space $\M_{2}(NDp^{\infty},L)$.\\

Now let $f$ be the new form of level $N$ in the introduction. Recall
the linear form $l_{f}$ defined on the space $\M_{2}(Np^{\infty},L)$
that interpolates "Peterson inner product with $f$". We would like
to take a trace to pass from forms of level $NDp^{\infty}$ to forms
of level $Np^{\infty}$.\\

For any integer $r>0$, the trace map
$\tr_{D}:\M_{2}(NDp^{r})\rightarrow\M_{2}(Np^{r})$ is defined in the
usual way, i.e. for any form $\phi\in\M_{2}(NDp^{r})$ we have
$\tr_{D}(\phi)(g)=\sum\limits_{h\in K_{0}(NDp^{r})\backslash
K_{0}(Np^{r})}\phi(gh)$. But under our hypothesis that $N$, $D$, $p$
are prime to each other, the set $K_{0}(NDp^{r})\backslash
K_{0}(Np^{r})$ can be identified with $K_{0}(D)\backslash K_{0}(1)$,
which does not depend on $r$, hence the trace maps $\tr_{D}$ on
different levels are coherent, and induce a map, which we still
denote by $\tr_{D}$, from the space $\M_{2}(NDp^{\infty},L)$ to the
space $\M_{2}(Np^{\infty},L)$.\\

Denote by $\mathbf{K}$ the pseudo-measure
$\tr_{D}(\widetilde{\mathbf{K}})$ and by $\mathbf{K}_{f}$ the
pseudo-measure $l_{f}(\mathbf{K})$. That is to say,
\begin{eqnarray*}
\mathbf{K}(\phi) &:=& \tr_{D}(\widetilde{\mathbf{K}}(\phi)),\\
\mathbf{K}_{f}(\phi) &:=& l_{f}(\mathbf{K}(\phi)).
\end{eqnarray*}
The pseudo-measure $\mathbf{K}$ is called the \textbf{analytic
kernel}.\\

We have the following interpolation property:

\begin{prop}
For any finite order character $\chi$ on $G$ of conductor $C\mid
p^\infty$, we have:
$$
\int\limits_{G}\chi\mathbf{K}_{f}=\frac{r(\chi)D_F^2\left|D\N(C)
\right|_\infty^{1/2}\chi^{-1}(\dEF)}{(8\pi^2)^n(f,f)}\left(\prod
\limits_{\PP_i\mid p}\frac{\alpha_i^{2-\val_{\PP_i}(C)}\prod
\limits_{q\mid\PP_i}(\alpha_i-\chi^{-1}(q))}{(\alpha_i^2-1)(\alpha_i^2-
\left|\PP_i\right|_\infty)}\right)\cdot L(1,f,\chi^{-1}).
$$
\end{prop}

This proposition is proved by Disegni in \cite{Disegni13} section
4.4, using his construction of Eisenstein series. But since the two
definitions lead to the same form, and the interpolation formula is
not needed in the following work, we won't reprove it here.\\

Now that we have the pseudo-measure $\mathbf{K}_{f}$, we may define
the $p$-adic $L$ function associated to $f$ and $E$ as the Mellin
transform of $\mathbf{K}_{f}$. More precisely, let $\mathfrak{X}$ be
the rigid analytic space $\Hom_{\mathrm{cont}}(G,\Cp^{\times})$,
then the pseudo-measure $\mathbf{K}_{f}$ induces a meromorphic
function $L_{p}$ on $\mathfrak{X}$:
$$
L_{p}(x):=\int\limits_{G}x\mathbf{K}_{f}.
$$
The only possible poles of $L_{p}$ are on the line
$\left\{x\in\mathfrak{X}:x_{F}=\xi_{p}\right\}$, which come from the
pole of the Deligne-Ribet $L$ function. In particular, the $L$
function is analytic on the anti-cyclotomic line
$\left\{x\in\mathfrak{X}:x=\overline{x}\right\}$. Here the bar means
action of the non-trivial element of $\Gal(E/F)$.\\

Suppose that $x_{0}\in\mathfrak{X}$ is not a pole of $L_{p}$, and
let $x$ be an element in the "neutral component" of $\mathfrak{X}$
(i.e. $x$ is trivial on the torsion subgroup of $G$). The derivative
of $L_{p}$ on the point $x_{0}$ along the direction $x$ is by
definition:
\begin{eqnarray*}
L_{p,x}'(x_{0}) &:=& \frac{d}{ds}L_{p}(x_{0}x^{s})\big|_{s=0}\\
&=&
\lim\limits_{s\rightarrow0}\frac{L_{p}(x_{0}x^{s})-L_{p}(x_{0})}{s},
\end{eqnarray*}
here the character $x^{s}$ is understood as $\exp(s\log x)$.\\

It is then easy to calculate that:
\begin{eqnarray*}
L_{p,x}'(x_{0}) &=&
\frac{d}{ds}\int\limits_{G}x_{0}x^{s}\mathbf{K}_{f}\big|_{s=0}\\
&=& \int\limits_{G}x_{0}x^{s}\log x\mathbf{K}_{f}\big|_{s=0}\\
&=& \int\limits_{G}x_{0}\log x\mathbf{K}_{f}.
\end{eqnarray*}

\section{$q$-expansion of the analytic kernel}

Since the pseudo-measure $\mathbf{K}$ takes values in the space of
modular forms, for each ideal $I$ of $\OF$ there are pseudo-measures
$a_{0}(I,\mathbf{K})$ and $a(I,\mathbf{K})$ taking values in $L$,
which correspond to the Fourier coefficients of $\mathbf{K}$. In
this section we are going to give explicit formulas for these
coefficients.\\

More precisely, let $\chi$ be a character mod $p^{r}$, we are going
to calculate the Fourier coefficients of the form
$\mathbf{K}(\chi)$. By definition, this form is:
$\tr_{D}\left[\Theta(\chi)\V_{N}
E_{1/2}(\epsilon\chi_{F}^{-1})\right]$. Note that a complete set of
representatives of $K(1)/K(D)$ is given by:
$$
K(1)/K(D)=\bigsqcup_{D_{1}\mid
D}\left\{\gamma_{D_{1},a}:a\in\widehat{\OF}/D_{1}\right\},
$$
where $\gamma_{D_{1},a}$ is the matrix whose $v$-th component is $
\begin{pmatrix}
a & 1\\
1 &
\end{pmatrix}$ if $v$ divides $D_{1}$ and is the identity matrix if
$v$ does not divide $D_{1}$.\\

Fix $d_{1}\in\Afx$ a generator of $D_{1}$. Denote by
$\gamma_{D_{1}}$ the matrix whose $v$-th component is $
\begin{pmatrix}
 & 1\\
d_{1} &
\end{pmatrix}$ if $v$ divides $D_{1}$ and is the identity matrix if
$v$ does not divide $D_{1}$. Note the identity:
$d_{1}\gamma_{D_{1},a}=
\begin{pmatrix}
d_{1} & a\\
 & 1
\end{pmatrix}\gamma_{D_{1}}$. Write $X$ for the form
$\Theta(\chi)\V_{N}E_{1/2}(\epsilon\chi_{F}^{-1})$. Since the form
$X$ has trivial central character, we have:
\begin{eqnarray*}
\mathbf{K}(\chi)(g) &=& \sum\limits_{D_{1}\mid
D}\sum\limits_{a\in\widehat{\OF}/D_{1}}X(g\gamma_{D_{1},a})\\
&=& \sum\limits_{D_{1}\mid
D}\sum\limits_{a\in\widehat{\OF}/D_{1}}X\left(g
\begin{pmatrix}
d_{1} & a\\
 & 1
\end{pmatrix}\gamma_{D_{1}}\right)\\
&=& \sum\limits_{D_{1}\mid
D}\sum\limits_{a\in\hat{O}/D_{1}}X_{D_{1}}\left(g
\begin{pmatrix}
d_{1} & a\\
 & 1
\end{pmatrix}\right)\\
&=& \sum\limits_{D_{1}\mid D}Y_{D_{1}}(g),
\end{eqnarray*}
here we have defined the forms $X_{D_{1}}$ and $Y_{D_{1}}$ by:
\begin{eqnarray*}
X_{D_{1}}(g) &:=& X(g\gamma_{D_{1}}),\\
Y_{D_{1}}(g) &:=&
\sum\limits_{a\in\widehat{\OF}/D_{1}}X_{D_{1}}\left(g
\begin{pmatrix}
d_{1} & a\\
 & 1
\end{pmatrix}\right).
\end{eqnarray*}

To calculate the Fourier coefficients of the form $\mathbf{K}(\chi)$
(i.e. the Whittaker functions $W\left(
\begin{pmatrix}
y & \\
 & 1
\end{pmatrix},\mathbf{K}(\chi)\right)$), it then suffices to
calculate the Whittaker functions of each $Y_{D_{1}}$. But it is
easy to see that:
$$
W\left(
\begin{pmatrix}
y & \\
 & 1
\end{pmatrix},Y_{D_{1}}\right)=
\begin{cases}
\left|D_{1}\right|_{\infty}W\left(
\begin{pmatrix}
d_{1}y & \\
 & 1
\end{pmatrix},X_{D_{1}}\right), & \textrm{if }yd_{F}\textrm{ is
integral};\\
0, & \textrm{else}.
\end{cases}
$$
Thus it suffices to calculate the Whittaker functions of
$X_{D_{1}}$.\\

Since $X_{D_{1}}$ is the product of the two forms:
\begin{eqnarray*}
\Theta_{D_{1}}(g) &:=& \Theta(\chi)(g\gamma_{D_{1}}),\\
E_{D_{1}}(g) &:=&
\V_{N}E_{1/2}(\epsilon\chi_{F}^{-1})(g\gamma_{D_{1}}),
\end{eqnarray*}
we see that the Whittaker function of $X_{D_{1}}$ can be obtained
from the Whittaker functions of $\Theta_{D_{1}}$ and $E_{D_{1}}$. We
now do the calculations one by one.

\subsection{The Whittaker function of $\Theta_{D_{1}}$}

The theta functions are defined in section \ref{Thetasection}. By
definition, we have:
$$
\Theta_{D_{1}}(g)=\sum\limits_{\sigma\in\AEx/\Ex(1+p^{r}
\widehat{\OEE})^{\times}E_{\infty}^{\times}}\chi(\sigma)^{-1}
\theta_{\sigma}(g\gamma_{D_{1}}),
$$
where:
\begin{eqnarray*}
\theta_{\sigma}(g\gamma_{D_{1}}) &=& \sum\limits_{(t,u)\in
U\backslash(E\times\Fx)}g\gamma_{D_{1}}\sigma\rho(t,u)\\
&=& \sum\limits_{(t,u)\in
U\backslash(E\times\Fx)}g\sigma\gamma_{D_{1}}\rho(t,u).
\end{eqnarray*}
Here we used the fact that the actions of $\GL(\A)$ and $\AEx$
commute with each other.\\

Writing $W_{\sigma}$ for the Whittaker function
$W(\cdot,\theta_{\sigma}(g\gamma_{D_{1}}))$, we have:
\begin{eqnarray*}
W_{\sigma}
\begin{pmatrix}
y & \\
 & 1
\end{pmatrix} &=&
\int\limits_{\A/F}\sum\limits_{U\backslash(E\times\Fx)}\left[
\begin{pmatrix}
1 & x\\
 & 1
\end{pmatrix}
\begin{pmatrix}
y & \\
 & 1
\end{pmatrix}\sigma\gamma_{D_{1}}\rho\right](t,u)\psi(-x)dx\\
&=& \sum\limits_{U\backslash(E\times\Fx)}\left[
\begin{pmatrix}
y & \\
 & 1
\end{pmatrix}\sigma\gamma_{D_{1}}\rho\right](t,u)\int
\limits_{\A/F}\psi(x(u\N(t)-1))dx\\
&=& \sum\limits_{t\in U\backslash\Ex}\left[
\begin{pmatrix}
y & \\
 & 1
\end{pmatrix}\sigma\gamma_{D_{1}}\rho\right](t,\N(t)^{-1}).
\end{eqnarray*}

The function $
\begin{pmatrix}
y & \\
 & 1
\end{pmatrix}\sigma\gamma_{D_{1}}\rho$
decomposes into product of local terms:
$$
\left[
\begin{pmatrix}
y & \\
 & 1
\end{pmatrix}\sigma\gamma_{D_{1}}\rho\right](t,\N(t)^{-1})=
\prod\limits_{v}\left[
\begin{pmatrix}
y_{v} & \\
 & 1
\end{pmatrix}\sigma_{v}(\gamma_{D_{1}})_{v}\rho_{v}\right]
(t,\N(t)^{-1})=:\prod\limits_{v}B_{v}.
$$
Here the product ranges over all places $v$ of $F$. We now calculate
each local term $B_{v}$.

\paragraph{Case 1: $v\mid\infty$}

Here both $\sigma_{v}$ and $(\gamma_{D_{1}})_{v}$ act trivially, so
we have:
$$
B_{v}=\mathds{1}_{y_{v}>0}\left|y\right|_{v}^{1/2}e^{-2\pi y}.
$$

\paragraph{Case 2: $v\nmid D_{1}\infty$}

Here again $(\gamma_{D_{1}})_{v}$ acts trivially, so we have:
$$
B_{v}=\mathds{1}_{\N(\sigma t)^{-1}\in y\dF\Ovx}\mathds{1}_{(\sigma
t)^{-1}\in(1+p^{r}\OEv)}\left|y\right|_{v}^{1/2}.
$$

\paragraph{Case 3: $v\mid D_{1}$}

In this case we have $(\gamma_{D_{1}})_{v}=
\begin{pmatrix}
 & 1\\
d_{1} &
\end{pmatrix}=
\begin{pmatrix}
 & 1\\
-1 &
\end{pmatrix}
\begin{pmatrix}
-d_{1} & \\
 & 1
\end{pmatrix}$, so that (notice that $1+p^{r}\OEv=\OEv$):
\begin{eqnarray*}
(\gamma_{D_{1}})_{v}\rho(t,u) &=&
\gamma_{v}\left|u\right|_{v}\left|D_{1}\right|_{v}^{1/2}
\int\limits_{E_{v}}\mathds{1}_{su\in\OEv}\mathds{1}_{d_{1}^{-1}u
\in\dF\Ovx}\psi(u\N(t,s))ds\\
&=& \gamma_{v}\left|u\right|_{v}\left|D_{1}\right|_{v}^{1/2}
\mathds{1}_{u\in
d_{1}\dF\Ovx}\int\limits_{u^{-1}\OEv}\psi(u\N(t,s))ds,
\end{eqnarray*}
hence:
$$
B_{v}=\left|y\right|_{v}^{1/2}\epsilon_{v}(y)\gamma_{v}
\left|(y\N(\sigma t))^{-1}\right|_{v}\left|D_{1}\right|_{v}^{1/2}
\mathds{1}_{(y\N(\sigma t))^{-1}\in
d_{1}\dF\Ovx}\int\limits_{y\N(\sigma t)\OEv}\psi((y\N(\sigma
t))^{-1}\N(y\sigma t,s))ds.
$$
Notice that:
\begin{eqnarray*}
(y\N(\sigma t))^{-1}\N(y\sigma t,s) &=& \N(\sigma
t)^{-1}(\overline{\sigma t}s+\sigma t\overline{s})\\
&=& \tr_{E/F}((\sigma t)^{-1}s),
\end{eqnarray*}
and since $\psi(\tr((\sigma t)^{-1}s))=\psi_{E}((\sigma t)^{-1}s)$,
we get:
\begin{eqnarray*}
B_{v} &=& \epsilon_{v}(y)\gamma_{v}\left|\N(\sigma
t)\right|_{v}^{-1}\left|D_{1}\right|_{v}^{1/2}\left|y
\right|_{v}^{-1/2}\mathds{1}_{\N(\sigma t)^{-1}\in y\dF
d_{1}\Ovx}\int\limits_{y\N(\sigma t)\OEv}\psi_{E}((\sigma
t)^{-1}s)ds\\
&=& \epsilon_{v}(y)\gamma_{v}\left|\N(\sigma
t)\right|_{v}^{-1}\left|D_{1}\right|_{v}^{1/2}\left|y
\right|_{v}^{-1/2}\mathds{1}_{\left|\N(\sigma t)\right|_{v}^{-1}=
\left|y\dF D_{1}\right|_{v}}\mathds{1}_{(\sigma t)^{-1}y\N(\sigma
t)\in d_{E}^{-1}\OEv}\left|y\N(\sigma
t)\right|_{E,v}\left|d_{E}\right|_{E,v}^{1/2}\\
&=& \epsilon_{v}(y)\gamma_{v}\left|y\right|_{v}^{1/2}
\mathds{1}_{\left|\N(\sigma t)\right|_{v}^{-1}= \left|y\dF
D_{1}\right|_{v}}\mathds{1}_{\left|\N(\sigma
t)\right|_{v}\leq\left|D_{1}\right|^{-1}\left|y\dF\right|^{-2}}\\
&=& \epsilon_{v}(y)\gamma_{v}\left|y\right|_{v}^{1/2}
\mathds{1}_{\left|\N(\sigma t)D_{1}\right|_{v}^{-1}=\left|y\dF
\right|_{v}}\mathds{1}_{\left|y\dF\right|_{v}\leq1}.
\end{eqnarray*}
Here we used the formula $d_{E,v}=d_{E/F,v}d_{F,v}$ and the fact
that $d_{E/F,v}$ is the maximal ideal of $\OEv$ in this case.

\paragraph{Putting together}

Denote by $\Delta_{1}$ the ideal of $\OEE$ such that
$\Delta_{1}^{2}=D_{1}\OEE$, and by $\delta_{1}\in\AEx$ a fixed
generator of $\Delta_{1}$. Combining the above results, we get:
\begin{eqnarray*}
&& \left[
\begin{pmatrix}
y & \\
 & 1
\end{pmatrix}\sigma\gamma_{D_{1}}\rho\right](t,\N(t)^{-1})\\
&=& \mathds{1}_{y_{\infty}>0}\left|y\right|_{\infty}^{1/2}e^{-2\pi
y_{\infty}}\prod\limits_{v\nmid
D_{1}}\left|y\right|_{v}^{1/2}\mathds{1}_{(\sigma
t)^{-1}\in(1+p^{r}\OEv)}\mathds{1}_{\N(\sigma t)^{-1}\in
y\dF\Ovx}\\
&& \prod\limits_{v\mid
D_{1}}\epsilon_{v}(y)\gamma_{v}\left|y\right|_{v}^{1/2}
\mathds{1}_{y\dF\in\Ov}\mathds{1}_{\N(\sigma t)^{-1}\in y\dF
d_{1}\Ovx}\\
&=& \left|y\right|^{1/2}e^{-2\pi
y_{\infty}}\mathds{1}_{y_{\infty}>0}\mathds{1}_{\N(\sigma
t\delta_{1})^{-1}\in y\dF\widehat{\OF}^{\times}}\mathds{1}_{(\sigma
t\delta_{1})^{-1}\in(1+p^{r}\widehat{\OEE})}\prod\limits_{v\mid
D_{1}}\epsilon_{v}(y)\gamma_{v}.
\end{eqnarray*}
Taking the sum over all $t$:
\begin{eqnarray*}
W_{\sigma}
\begin{pmatrix}
y & \\
 & 1
\end{pmatrix} &=&
\mathds{1}_{y_{\infty}>0}\left|y\right|^{1/2}e^{-2\pi
y_{\infty}}\prod\limits_{v\mid D_{1}}\epsilon_{v}(y)\gamma_{v}\\
&& \sum\limits_{t\in U\backslash\Ex}\mathds{1}_{\N(\sigma
t\delta_{1})^{-1}\in y\dF\widehat{\OF}^{\times}}\mathds{1}_{(\sigma
t\delta_{1})^{-1}\in(1+p^{r}\widehat{\OEE})}\\
&=& \mathds{1}_{y_{\infty}>0}\left|y\right|^{1/2}e^{-2\pi
y_{\infty}}\prod\limits_{v\mid
D_{1}}\epsilon_{v}(y)\gamma_{v}\sum\limits_{\substack{J\subseteq
\OEE\\\N(J)=y\dF}}\mathds{1}_{(\sigma\delta_{1})^{-1}}(J),
\end{eqnarray*}
here $\mathds{1}_{(\sigma\delta_{1})^{-1}}$ is the characteristic
function of the class of $(\sigma\delta_{1})^{-1}$ in
$\AEx/\Ex(1+p^{r}\widehat{\OEE})^{\times}E_{\infty}^{\times}$, i.e.
it is equal to $1$ if and only if $J$ is prime to $p$ and a
generator $j\in\AEx$ of $J$ is in the same class as
$(\sigma\delta_{1})^{-1}$. Note that we used the bijection $t\mapsto
J=(\sigma t\delta_{1})^{-1}\OEE$ above.\\

Finally, taking the sum over all $\sigma$:
\begin{eqnarray*}
W\left(
\begin{pmatrix}
y & \\
 & 1
\end{pmatrix},\Theta_{D_{1}}\right) &=&
\sum\limits_{\sigma}\chi(\sigma)^{-1}W_{\sigma}
\begin{pmatrix}
y & \\
 & 1
\end{pmatrix}\\
&=& \mathds{1}_{y_{\infty}>0}\left|y\right|^{1/2}e^{-2\pi
y_{\infty}}\prod\limits_{v\mid
D_{1}}\epsilon_{v}(y)\gamma_{v}\sum\limits_{\substack{J\subseteq
\OEE\\\N(J)=y\dF}}
\sum\limits_{\sigma}\mathds{1}_{(\sigma\delta_{1})^{-1}}(J)
\chi(\sigma)^{-1}\\
&=& \mathds{1}_{y_{\infty}>0}\left|y\right|^{1/2}e^{-2\pi
y_{\infty}}\chi_{[p]}(\Delta_{1})\prod\limits_{v\mid
D_{1}}\epsilon_{v}(y)\gamma_{v}\sum\limits_{\substack{J
\subseteq\OEE\\\N(J)=y\dF}}\chi_{[p]}(J).
\end{eqnarray*}

It is easy to see that the constant term $W_{0}\left(
\begin{pmatrix}
y & \\
 & 1
\end{pmatrix},\Theta_{D_{1}}\right)$ vanishes, since a similar
calculation shows that the constant term of each
$\theta_{\sigma}(g\gamma_{D_{1}})$ vanishes.

\subsection{The Whittaker function of $E_{D_{1}}$}

Write $\eta$ for the character $\epsilon\chi_{F}^{-1}$ and $C$ for
the conductor of $\chi_{F}$ (hence the conductor of $\eta$ is $CD$).
We begin with a complex variable $s$ with $\mathrm{Re}(s)>1$ and
write $\mathcal{E}_{s}$ for the form
$\V_{N}E_{s}(\eta)(g\gamma_{D_{1}})$ (thus $E_{D_{1}}$ is equal to
$\mathcal{E}_{1/2}$). By definition, we have:
$$
\mathcal{E}_{s}=2^{-n}L_{CD}(1-2s,\eta)\sum \limits_{\gamma\in
B(F)\backslash\GL(F)}Q(\gamma g),
$$
where the function $Q=\prod\limits_{v}Q_{v}$ satisfies $Q\left(
\begin{pmatrix}
a & b\\
 & d
\end{pmatrix}k\right)=\left|\frac{a}{d}\right|^{s}\eta(d)Q(k)$ for
any $k$ in the maximal compact subgroup, and each $Q_{v}$ is defined
by:
\begin{itemize}
\item if $v\mid\infty$, then we have
$Q_{v}(r(\theta))=e^{i\theta}$;
\item if $v\nmid pD\infty$, then for $k=
\begin{pmatrix}
u & v\\
w & t
\end{pmatrix}\in\GL(\Ov)$ we have $Q_{v}\left(k
\begin{pmatrix}
n_{f} & \\
 & 1
\end{pmatrix}\right)=\left|N\right|_{v}^{1/2}$, here
$n_{f}\in\Afx$ denotes a generator of $N$;
\item if $v\mid CD$ but $v\nmid D_{1}$, then for $k=
\begin{pmatrix}
u & v\\
w & t
\end{pmatrix}\in\GL(\Ov)$ we have
$$
Q_{v}(k)=
\begin{cases}
\eta_{v}(t), & \textrm{if }(CD)_{v}\mid w;\\
0, & \textrm{else};
\end{cases}
$$
\item if $v\mid D_{1}$, then for $k=
\begin{pmatrix}
u & v\\
w & t
\end{pmatrix}\in\GL(\Ov)$ we have
$$
Q_{v}(k)=
\begin{cases}
0, & \textrm{if }(D_{1})_{v}\mid w;\\
\left|D_{1}\right|_{v}^{s}\eta_{v}(w), & \textrm{else};
\end{cases}
$$
\item if $v\mid p$ but $v\nmid C$, then for $k=
\begin{pmatrix}
u & v\\
w & t
\end{pmatrix}\in\GL(\Ov)$ we have
$$
Q_{v}(k)=
\begin{cases}
1-\left|\pi_{v}\right|, & \textrm{if }\left|w\right|_{v}=1;\\
1-\eta(\pi_{v})\left|\pi_{v}\right|^{1-2s}, & \textrm{else}.
\end{cases}
$$
\end{itemize}

The following calculation of Fourier coefficients is standard.
\begin{eqnarray*}
W\left(
\begin{pmatrix}
y & \\
 & 1
\end{pmatrix},\mathcal{E}_{s}\right) &=&
2^{-n}L_{CD}(1-2s,\eta)\int\limits_{\A/F}\left(Q
\begin{pmatrix}
y & x\\
 & 1
\end{pmatrix}+\sum\limits_{u\in F}Q
\begin{pmatrix}
 & 1\\
y & u+x
\end{pmatrix}\right)\psi(-x)dx\\
&=& 2^{-n}L_{CD}(1-2s,\eta)\int\limits_{\A}Q
\begin{pmatrix}
 & 1\\
y & x
\end{pmatrix}\psi(-x)dx\\
&=& 2^{-n}L_{CD}(1-2s,\eta)\left|y\right|^{1-s}\eta(y)\prod
\limits_{v}V_{v}(y_{v}),
\end{eqnarray*}
where the functions $V_{v}$ are defined as:
$$
V_{v}(y):=\int\limits_{F_{v}}Q_{v}
\begin{pmatrix}
 & 1\\
1 & x
\end{pmatrix}\psi_{v}(-yx)dx.
$$

The constant term can be computed similarly:
$$
W_{0}\left(
\begin{pmatrix}
y & \\
 & 1
\end{pmatrix},\mathcal{E}_{s}\right)=2^{-n}L_{CD}(1-2s,\eta)
\left(\left|y\right|^{s}Q
\begin{pmatrix}
1 & \\
 & 1
\end{pmatrix}+\left|y\right|^{1-s}\eta(y)\prod\limits_{v}
V_{v}(0)\right).
$$
Note that the value of $Q
\begin{pmatrix}
1 & \\
 & 1
\end{pmatrix}$ is easily seen to be
$\mathds{1}_{D_{1}=1}\left|N\right|_{f}^{1/2-s}\prod
\limits_{\substack{v\mid p\\v\nmid
C}}(1-\eta_{v}(\pi_{v})\left|\pi_{v}\right|_{v}^{1-2s})$.

It then remains to compute each local function $V_{v}$.

\paragraph{Case 1: $v\mid\infty$}

In this case the standard calculation gives:
$$
V_{v}(y)=\int\limits_{\Fv}\frac{e^{-2\pi
iyx}}{(x^{2}+1)^{s-1/2}(x+i)}dx.
$$
If $y$ is non-zero, then this function can be analytically continued
to all $s\in\C$, and decreases exponentially with respect to
$\left|y\right|$. Furthermore, when $s=1/2$, we have:
$$
V_{v}(y)=
\begin{cases}
(-2\pi i)e^{-2\pi y}, & \textrm{if }y>0;\\
0, & \textrm{if }y<0.
\end{cases}
$$
If y is zero, then we have:
$$
V_{s}(0)=-i\pi^{1/2}\Gamma(s)\Gamma(s+1/2)^{-1}.
$$

\paragraph{Case 2: $v\nmid pD\infty$}

In this case we have:
$$
Q_{v}
\begin{pmatrix}
 & 1\\
1 & x
\end{pmatrix}=
\begin{cases}
\eta_{v}(N)^{-1}\left|N\right|_{v}^{s+1/2}, & \textrm{if
}\left|Nx\right|\leq1;\\
\eta_{v}(x)\left|N\right|_{v}^{1/2-s}\left|x\right|_{v}^{-2s}, &
\textrm{else}.
\end{cases}
$$
It follows that:
\begin{eqnarray*}
V_{v}(y) &=& \eta_{v}(N)^{-1}\left|N\right|_{v}^{s+1/2}
\int\limits_{\left|x\right|_{v}\leq\left|N\right|_{v}^{-1}}
\psi_{v}(-yx)dx\\
&& +\left|N\right|_{v}^{1/2-s}\sum\limits_{i>\val_{v}(N)}
\int\limits_{\left|x\right|_{v}=\left|\pi_{v}\right|_{v}^{-i}}
\eta_{v}(x)\left|x\right|_{v}^{-2s}\psi_{v}(-yx)dx\\
&=& \eta_{v}(N)^{-1}\left|N\right|_{v}^{s+1/2}S_{\val_{v}(N)}+
\left|N\right|_{v}^{1/2-s}\sum\limits_{i>\val_{v}(N)}\eta_{v}
(\pi_{v})^{-i}\left|\pi_{v}\right|_{v}^{2is}(S_{i}-S_{i-1}),
\end{eqnarray*}
where $S_{i}$ is defined as:
$$
S_{i}:=\int\limits_{\left|x\right|\leq\left|\pi_{v}\right|_{v}^{-i}}
\psi_{v}(-yx)dx=
\begin{cases}
\left|\pi_{v}\right|_{v}^{-i}\left|\dF\right|_{v}^{1/2}, &
\textrm{if }\left|y\dF\right|_{v}\leq\left|\pi_{v}\right|_{v}^{i};\\
0, & \textrm{else}.
\end{cases}
$$
Hence the function $V_{v}$ is non-zero only if
$\left|y\dF\right|_{v}\leq\left|N\right|_{v}$, and in this case we
have (making the convention that $S_{i}=0$ for $i<\val_{v}(N)$):
\begin{eqnarray*}
V_{v}(y) &=&
\left|N\right|_{v}^{1/2-s}\sum\limits_{i}\eta_{v}(\pi_{v})^{-i}\left|
\pi_{v}\right|_{v}^{2is}(S_{i}-S_{i-1})\\
&=& \left|N\right|_{v}^{1/2-s}(1-\eta_{v}(\pi_{v})^{-1}\left|\pi_{v}
\right|_{v}^{2s})\sum\limits_{i}\eta_{v}(\pi_{v})^{-i}\left|\pi_{v}
\right|_{v}^{2is}
S_{i}\\
&=& \left|N\right|_{v}^{1/2-s}(1-\eta_{v}(\pi_{v})^{-1}\left|\pi_{v}
\right|_{v}^{2s})\left|\dF\right|_{v}^{1/2}\sum\limits_{i=
\val_{v}(N)}^{\val_{v}(y\dF)}(\eta_{v}(\pi_{v})^{-1}\left|
\pi_{v}\right|_{v}^{2s-1})^{i}.
\end{eqnarray*}

\paragraph{Case 3: $v\mid CD$ but $v\nmid D_{1}$}

In this case a similar calculation shows that the function $V_{v}$
is reduced to a Gauss sum: it is zero if $y\dF$ is not integral or
if $y=0$; otherwise we have:
$$
V_{v}(y)=\left|CD\right|_{v}^{2s-1/2}\left|y\dF\right|_{v}^{2s-1}
\left|\dF\right|_{v}^{1/2}\eta_{v}(-y)^{-1}r_{v}(\eta,\psi),
$$
where $r_{v}(\eta,\psi)$ is the root number:
$$
r_{v}(\eta,\psi)=\left|CD\dF\right|_{v}^{1/2}\int
\limits_{\left|x\right|_{v}=\left|CD\dF\right|_{v}^{-1}}\eta_{v}(x)
\psi_{v}(x)dx.
$$
Note that for places $v\mid D$ we have
$r_{v}(\eta,\psi)=r_{v}(\epsilon,\psi)\chi_{F,v}(D\dF)$, and
$r_{v}(\epsilon,\psi)$ is equal to $\gamma_{v}$, the Weil index.

\paragraph{Case 4: $v\mid D_{1}$}

In this case we have:
$$
Q_{v}
\begin{pmatrix}
 & 1\\
1 & x
\end{pmatrix}=
\begin{cases}
\left|D_{1}\right|_{v}^{s}, & \textrm{if }\left|x\right|_{v}\leq1;\\
0, & \textrm{else},
\end{cases}
$$
so that:
$$
V_{v}(y)=\int\limits_{\Ov}\left|D_{1}\right|_{v}^{s}\psi_{v}(-yx)dx=
\begin{cases}
\left|D_{1}\right|_{v}^{s}\left|\dF\right|_{v}^{1/2}, & \textrm{if
}\left|y\dF\right|_{v}\leq1;\\
0, & \textrm{else}.
\end{cases}
$$

\paragraph{Case 5: $v\mid p$ but $v\nmid C$}

In this case we have:
$$
Q_{v}
\begin{pmatrix}
 & 1\\
1 & x
\end{pmatrix}=
\begin{cases}
1-\left|\pi_{v}\right|_{v}, & \textrm{if }\left|x\right|_{v}\leq1;\\
\eta_{v}(x)\left|x\right|_{v}^{-2s}(1-\eta_{v}(\pi_{v})\left|\pi_{v}
\right|_{v}^{1-2s}), & \textrm{else}.
\end{cases}
$$
From this we see that $V_{v}(y)$ is non-zero only if $y\dF$ is
integral, in which case we have:
$$
V_{v}(y)=(1-\left|\pi_{v}\right|_{v})\left|\dF\right|_{v}^{1/2}+
(1-\eta_{v}(\pi_{v})\left|\pi_{v}\right|_{v}^{1-2s})\sum\limits_{i>0}
\eta_{v}(\pi_{v})^{-i}\left|\pi_{v}\right|_{v}^{2is}
(S_{i}-S_{i-1}),
$$
where:
$$
S_{i}:=\int\limits_{\left|x\right|_{v}\leq\left|\pi_{v}
\right|_{v}^{-i}}\psi_{v}(-yx)=
\begin{cases}
\left|\pi_{v}\right|_{v}^{-i}\left|\dF\right|_{v}^{1/2}, &
\textrm{if }i\leq\val_{v}(y\dF);\\
0, & \textrm{else}.
\end{cases}
$$
After simplification, we get (when $y\dF$ is integral):
$$
V_{v}(y)=
\begin{cases}
(1-\eta_{v}(\pi_{v})^{-1}\left|\pi_{v}\right|_{v}^{2s})\left|\dF
\right|_{v}^{1/2}\eta_{v}(y\dF)^{-1}\left|y\dF
\right|_{v}^{2s-1}, & \textrm{if }y\neq0;\\
0, & \textrm{if }y=0.
\end{cases}
$$

\paragraph{Putting together}

Combining the above results, we may now compute the Whittaker
function of $\mathcal{E}_{s}$:
\begin{eqnarray*}
&& W\left(
\begin{pmatrix}
y & \\
 & 1
\end{pmatrix},\mathcal{E}_{s}\right)\\
&=& 2^{-n}L_{CD}(1-2s,\eta)
\left|y\right|^{1-s}\eta(y)\prod\limits_{v\mid\infty}V_{v}(y_{v})\\
&& \prod \limits_{v\nmid
pD\infty}\mathds{1}_{\left|y\dF\right|_{v}\leq\left|N
\right|_{v}}\left|N\right|_{v}^{1/2-s}(1-\eta_{v}(\pi_{v})^{-1}\left|
\pi_{v}\right|_{v}^{2s})\left|\dF\right|_{v}^{1/2}
\sum\limits_{i=\val_{v}(N)}^{\val_{v}(y\dF)}(\eta_{v}
(\pi_{v})^{-1}\left|\pi_{v}\right|_{v}^{2s-1})^{i}\\
&& \prod\limits_{\substack{v\mid CD\\v\nmid
D_{1}}}\mathds{1}_{\left|y\dF\right|_{v}\leq1}\left|CD
\right|_{v}^{2s-1/2}\left|y\dF\right|_{v}^{2s-1}\left|\dF
\right|_{v}^{1/2}\eta_{v}(-y)^{-1}r_{v}(\eta,\psi)\\
&& \prod\limits_{v\mid
D_{1}}\mathds{1}_{\left|y\dF\right|_{v}\leq1}\left|D_{1}
\right|_{v}^{s}\left|\dF\right|_{v}^{1/2}\prod\limits_{
\substack{v\mid p\\v\nmid C}}\mathds{1}_{\left|y\dF\right|_{v}\leq1}
(1-\eta_{v}(\pi_{v})^{-1}\left|\pi_{v}\right|_{v}^{2s})
\left|\dF\right|_{v}^{1/2}\eta_{v}(y\dF)^{-1}\left|y\dF
\right|_{v}^{2s-1}\\
&=& 2^{-n}L_{CD}(1-2s,\eta)L_{CD}(2s,\eta^{-1})^{-1}
\left|y\right|^{1-s}\eta(y)\mathds{1}_{N\mid
y_{f}\dF}\left|N\right|_{f}^{1/2-s}\left|D_{1}\right|_{f}^{1/2-s}
\left|CD\right|_{f}^{2s-1/2}\left|\dF\right|_{f}^{1/2}\\
&& \prod\limits_{v\mid\infty}V_{v}(y_{v})\prod \limits_{v\nmid
pD\infty}\sum\limits_{i=\val_{v}(N)}^{\val_{v}(y\dF)}(\eta_{v}
(\pi_{v})^{-1}\left|\pi_{v}\right|_{v}^{2s-1})^{i}\\
&&\prod \limits_{\substack{v\mid CD\\v\nmid
D_{1}}}\left|y\dF\right|_{v}^{2s-1}\eta_{v}(-y)^{-1}r_{v}
(\eta,\psi)\prod\limits_{\substack{v\mid p\\v\nmid C}}
\eta_{v}(y\dF)^{-1}\left|y\dF\right|_{v}^{2s-1}.\\
\end{eqnarray*}
Using the functional equation:
\begin{eqnarray*}
&& L_{CD}(1-2s,\eta)/L_{CD}(2s,\eta^{-1})\\
&=& \left|CD\dF
\right|_{\infty}^{2s-1/2}(-i\pi^{1/2-2s}\Gamma(1-s)^{-1}
\Gamma(s+1/2))^{n}\prod\limits_{v\nmid
CD\infty}\eta_{v}(\dF)\prod\limits_{v\mid CD}r_{v}(\eta,\psi)^{-1},
\end{eqnarray*}
and then specialize to $s=1/2$, we simplify the formula to:
\begin{eqnarray*}
&& W\left(
\begin{pmatrix}
y & \\
 & 1
\end{pmatrix},E_{D_{1}}\right)\\
&=& \mathds{1}_{y_{\infty}>0}\left|y\right|^{1/2}e^{-2\pi
y_{\infty}}\sum\limits_{\substack{J\subseteq\OF\\J\mid
y\dF/N}}\eta_{[pD]}(J)\prod\limits_{v\mid
D_{1}}\eta_{v}(-y)r_{v}(\eta,\psi)^{-1}\\
&=& \mathds{1}_{y_{\infty}>0}\left|y\right|^{1/2}e^{-2\pi
y_{\infty}}\sum\limits_{\substack{J\subseteq\OF\\J\mid
y\dF/N}}\chi^{-1}_{[p]}(J\OEE)\epsilon_{[D]}(J)\prod\limits_{v\mid
D_{1}}\epsilon_{v}(-y)r_{v}(\epsilon,\psi)^{-1}\chi_{F,v}^{-1}(y\dF
D_{1}).
\end{eqnarray*}

For the constant term, notice that for any place $v$ dividing $p$,
we always have $V_{v}(0)=0$. Hence:
\begin{eqnarray*}
W_{0}\left(
\begin{pmatrix}
y & \\
 & 1
\end{pmatrix},\mathcal{E}_{s}\right) &=&
2^{-n}L_{CD}(1-2s,\eta)\left|y\right|^{s}Q
\begin{pmatrix}
1 & \\
 & 1
\end{pmatrix}\\
&=& \mathds{1}_{D_{1}=1}2^{-n}\left|N\right|_{f}^{1/2-s}
\left|y\right|^{s}L_{pD}(1-2s,\eta).
\end{eqnarray*}
When specialized to $s=1/2$, this gives:
$$
W_{0}\left(
\begin{pmatrix}
y & \\
 & 1
\end{pmatrix},E_{D_{1}}\right)=\mathds{1}_{D_{1}=1}
\left|y\right|^{1/2}2^{-n}L_{pD}(0,\chi_{F}^{-1}\epsilon).
$$

\subsection{$q$-expansion of the analytic kernel}\label{analytic}

We are ready to calculate the $q$-expansion of the analytic kernel.
For each divisor $D_{1}$ of $D$ and each $y\in\Ax$ such that
$y_{\infty}>0$, we have:
\begin{eqnarray*}
&& W\left(
\begin{pmatrix}
y & \\
 & 1
\end{pmatrix},X_{D_{1}}\right)\\
&=&
\sum\limits_{\substack{\alpha+\beta=1\\\alpha,\beta\geq0}}W\left(
\begin{pmatrix}
\alpha y & \\
 & 1
\end{pmatrix},\Theta_{D_{1}}\right)W\left(
\begin{pmatrix}
\beta y & \\
 & 1
\end{pmatrix},E_{D_{1}}\right)\\
&=& \left|y\right|e^{-2\pi
y_{\infty}}\sum\limits_{\substack{\alpha+\beta=1\\
\alpha,\beta>0}}\sum\limits_{\substack{J\subseteq\OEE\\\N(J)=\alpha
y\dF\\(J,p)=1}}\sum\limits_{\substack{K\subseteq\OF\\K\mid\beta
y\dF/N\\(K,pD)=1}}\chi_{[p]}(J/\Delta_{1}K(\beta
y\dF)_{D{1}})\epsilon_{[D]}(K)\prod\limits_{v\mid
D_{1}}\epsilon_{v}(-\alpha\beta)\\
&& +\left|y\right|e^{-2\pi
y_{\infty}}\mathds{1}_{D_{1}=1}2^{-n}L_{pD}(0,\chi_{F}^{-1}
\epsilon)\sum\limits_{\substack{J\subseteq\OEE\\\N(J)=y\dF}}
\chi_{[p]}(J),
\end{eqnarray*}
here we write $(\beta y\dF)_{D_{1}}$ for the ideal
$\prod\limits_{v\mid D_{1}}(\beta y\dF)_{v}$ of $\OF$. The constant
term vanishes, because that
of $\Theta_{D_{1}}$ vanishes.\\

Finally, for any integral ideal $I$ of $\OF$, take $y\in\Ax$ such
that $y_{\infty}>0$ and $y\dF=I$, we may compute:
\begin{eqnarray*}
a(I,\mathbf{K}(\chi)) &=& \left|y\right|^{-1}e^{2\pi
y_{\infty}}W\left(
\begin{pmatrix}
y & \\
 & 1
\end{pmatrix},\mathbf{K}(\chi)\right)\\
&=& \left|y\right|^{-1}e^{2\pi y_{\infty}}\sum\limits_{D_{1}\mid
D}\left|D_{1}\right|_{\infty}W\left(
\begin{pmatrix}
d_{1}y & \\
 & 1
\end{pmatrix},X_{D_{1}}\right)\\
&=& \sum\limits_{D_{1}\mid
D}\sum\limits_{\substack{\alpha+\beta=1\\\alpha,\beta>0}}\sum
\limits_{\substack{J\subseteq\OEE\\\N(J)=\alpha
D_{1}I\\(J,p)=1}}\sum\limits_{\substack{K\subseteq\OF\\K\mid\beta
D_{1}I/N\\(K,pD)=1}}\chi_{[p]}(J/\Delta_{1}K(\beta
D_{1}I)_{D_{1}})\epsilon_{[D]}(K)\prod\limits_{v\mid
D_{1}}\epsilon_{v}(-\alpha\beta)\\
&& +2^{-n}L_{pD}(0,\chi_{F}^{-1}\epsilon)\sum\limits_{
\substack{J\subseteq\OEE\\\N(J)=I}}\chi_{[p]}(J).
\end{eqnarray*}
Since this formula is valid for all character $\chi$, by linearity
it is also valid for all functions on the group $G$.

\subsection{The central derivative}

Here we specialize to the case of the central derivative of the
$p$-adic $L$ function in the cyclotomic line. Writing $\xi_E$ for
the character $\xi_p\circ\N_{E/F}$ (the $p$-adic cyclotomic
character on $E$), we have:
$$
L'_{p,\xi_E}(1)=\int\limits_{G}(\log_{p}\circ\xi_E)\mathbf{K}_{f}
=l_f\left(\int\limits_{G}(\log_{p}\circ\xi_E)\mathbf{K}\right).
$$
We denote by $\Phi$ the $p$-adic modular form
$\int\limits_{G}(\log_{p}\circ\xi_{E})\mathbf{K}$. The calculations
in previous subsections apply and gives:
\begin{prop}\label{analyticcoeff}
If $I$ is an ideal of $\OF$ divisible by $p$, then the coefficient
$a(I,\Phi)$ of the $q$-expansion of the form $\Phi$ is given by
$a(I,\Phi)=\sum\limits_{v}a_v(I,\Phi)$, where the sum ranges over
finite places $v$, and each term $a_v(I,\Phi)$ is given by:
\begin{enumerate}
\item if $v$ is inert in $E$ and $q_v$ is the associated prime ideal, then:
$$
a_v(I,\Phi)=\sum\limits_{\substack{\alpha+\beta=1,
\alpha,\beta>0\\\epsilon_{w}(-\alpha\beta)=1,\forall w\mid
D\\N\mid\beta DI\\(p,\beta DI)=1}}\delta r(\alpha I)r(\beta
I/Nq_{v})\ord_{v}(\beta Iq_v/N)\log_p(\left|q_{v}\right|_\infty),
$$
here $\delta$ denotes $2^{\#\{v:v\mid D,v\mid\beta DI\}}$, and for
an ideal $M$ of $\OF$, $r(M)$ denotes the number of ideals of $\OEE$
with norm $M$;
\item if $v$ is ramified in $E$ and $q_v$ is the associated prime ideal,
then:
$$
a_v(I,\Phi)=\sum\limits_{\substack{\alpha+\beta=1,
\alpha,\beta>0\\\epsilon_{w}(-\alpha\beta)=1,\forall v\neq w\mid
D\\N\mid\beta DI\\\epsilon_{v}(-\alpha\beta)=-1\\(p,\beta
DI)=1}}\delta r(\alpha I)r(\beta I/N)\ord_{v}(\beta
Iq_v)\log_p(\left|q_{v}\right|_\infty);
$$
\item if $v$ is split in $E$, then $a_{v}(I,\Phi)=0$.
\end{enumerate}
\end{prop}

\begin{dem}
If we replace the character $\chi$ in the formula in subsection
\ref{analytic} by the function $\log_{p}\circ\xi_E$, then the sum
over all the $J$ becomes a multiplication by the factor $r(\alpha
I)$, because every $J$ has the same norm to $F$, and also note that
$r(\alpha D_{1}I)=r(\alpha I)$ for any $D_{1}\mid D$.\\

Since we only consider the case $p\mid I$, the contribution of the
second term is zero. Thus we have:
$$
a(I,\Phi)=\sum\limits_{\alpha+\beta=1}\sum\limits_{D_1}\sum\limits_{K}
r(\alpha I)\log_{p}(\xi_p(\alpha I/K^{2}(\beta
D_{1}I)_{D_1}^{2}))\epsilon_{[D]}(K)\prod\limits_{v\mid
D_{1}}\epsilon_{v}(-\alpha\beta).
$$
Now for any ideal $M$ of $\OF$, we have:
$$
\log_{p}(\xi_p(M))=\log_{p}\left(\prod\limits_{v}
\xi_p(q_{v}^{\ord_{v}(M)})\right)=
\sum\limits_{v}\ord_{v}(M)\log_{p}(\left|q_{v}\right|_\infty),
$$
thus the coefficient $a(I,\Phi)$ decomposes into a sum of local
terms $a_v(I,\Phi)$, and each local term can be computed.
\end{dem}

\newpage

\centerline{\Large{Part II. Geometric part}}

\section{Shimura curves and Heegner points}

In the following we assume that $\epsilon(N)=(-1)^{n-1}$. We are
going to introduce a Shimura curve $X$ and construct Heegner points
on the curve.

\subsection{Shimura curves}

Fix an infinite place $\tau$ of $F$. Under our hypothesis, there is
a quaternion algebra $B$ over $F$ which ramifies exactly at all
infinite places different from $\tau$ and all finite places $v$ such
that $\epsilon_v(N)=-1$.\\

By fixing an isomorphism $B_{\tau}\simeq M_{2}(\R)$, the group $\Bx$
acts on the left on the Poincar\'e double half plane
$\mathcal{H}^{\pm}:=\C-\R$ via the usual action of $\GL(\R)$. For
any open compact subgroup $K$ of
$\widehat{B}^{\times}/\widehat{F}^{\times}$, we then have a Shimura
curve $X_{K}$, whose complex points are given by:
$X_{K}(\C)=\Bx\backslash\mathcal{H}^{\pm}
\times\widehat{B}^{\times}/\widehat{F}^{\times}K$.

\begin{rem}
When $F=\Q$ there may be a finite set of cusps, but we will exclude
this case, which is the original work of Perrin-Riou.
\end{rem}

\subsection{Level structure and the Shimura curve $X$}\label{curveX}

By construction of the quaternion algebra $B$, the quadratic
extension $E$ is non-split on every ramified place of $B$, so we can
fix an embedding of $E$ into $B$ (which is unique up to conjugation
by $\Bx$), and thus view $E$ as a sub-algebra
of $B$.\\

In \cite{Zhang01} section 1.5, an order $R$ of $B$ of type $(N,E)$
is constructed, that is, the order $R$ contains $\OEE$ and has
discriminant $N$. Using this order $R$, we define the level
structure $K=\widehat{F}^{\times}\widehat{R}^{\times}$, and the
Shimura curve $X:=X_{K}$, whose complex points are given by:
$X(\C)=\Bx\backslash\mathcal{H}^{\pm}
\times\widehat{B}^{\times}/\widehat{F}^{\times}\widehat{R}^{\times}$.\\

The moduli interpretation problem of the curve $X$ is discussed in
\cite{Zhang01} section 1. There is a finite map from $X(\C)$ to
another Shimura curve $X'(\C)$ which parametrizes certain classes of
abelian varieties. This interpretation gives an integral model of
the curve $X$.

\subsection{CM points, Galois actions and Hecke
operators}\label{sectionCM}

Let $z_0$ be the unique point in the Poincar\'e upper half plane
which is fixed by the action of $\Ex$. A \textbf{CM point} on $X$ is
a point in $X(\C)$ which is represented by a pair
$(z_{0},b)\in\mathcal{H}^{\pm}\times\widehat{B}^{\times}$.\\

The Galois actions on CM points are described by Shimura's
reciprocity law, which states that all CM points are algebraic and
defined over the maximal abelian extension $E^{\mathrm{ab}}$ of $E$,
and for any element $a\in\widehat{E}^{\times}$, we have:
$$
\rec(a)(z_0,b)=(z_0,ab),
$$
here $\rec:\widehat{E}^\times\rightarrow\Gal(E^{\mathrm{ab}}/E)$ is
the Artin reciprocity map over $E$.\\

There are also Hecke operators defined on the group of divisors of
$X$, in a way similar to Heck operators on modular forms. Fix
isomorphisms $B_{v}\simeq M_{2}(F_{v})$ for split places $v$ of $B$.
Let $M$ be a non-zero ideal of $\OF$ prime to all rafimied places of
$B$, and let $U(M)$ be the same set as in subsection
\ref{heckeoperator}. Again write $U(M)$ as a disjoint union
$\bigsqcup\limits_{j}h_{j}K_{0}(N)$, and the Hecke operator $\TN(M)$
is defined as:
$$
\TN(M)(z,b)=\sum\limits_{j}(z,bh_j)=\sum\limits_{h\in
U(M)/K_{0}(N)}(z,bh).
$$
The Hecke operators are obviously multiplicative in $M$.\\

By Jacquet-Langlands theory, the Hecke algebra generated by these
Hecke operators is a quotient of the Hecke algebra generated by the
Hecke operators on modular forms of level $N$. Thus we may write
$\TN(I)$ unambiguously to represent both of the operators.

\subsection{Heegner divisors}\label{heegnerdivisor}

A CM point on the curve $X$ is called a "Heegner point" if it is of
conductor $1$, i.e. it is defined over the conductor $1$ ring class
field $H$ of $E$ (i.e. the abelian extension of $E$ such that
$\Gal(H/E)$ is isomorphic to
$\AEfx/\Ex\widehat{\OEE}^{\times}\widehat{F}^{\times}$ via class
field theory). The action of the group $\Gal(H/E)$ acts transitively
on the set of Heegner points.\\

More generally, if $C$ is a non-zero ideal of $\OF$, let $H[C]$ be
the ring class field of conductor $C$, i.e. the abelian extension of
$E$ such that
$\Gal(H[C]/E)\simeq\AEfx/\Ex(\widehat{\OF}+C\widehat{\OE})^{\times}
\widehat{F}^{\times}$. A CM point of conductor $C$ is a point that
is defined over the field $H[C]$.\\

In \cite{Zhang01}, a canonical divisor class of degree $1$,
$\xi\in\Pic(X)\otimes\Q$, is constructed (see section 4.1.4 and
Introduction of \textit{loc. cit.}),and one may use it to define a
map $\phi:X\rightarrow\Jac(X)\otimes\Q$, which sends any point $y$
to the class of $y-\xi$. The choice of a divisor representing the
Hodge class will be discussed later.\\

Now let $x$ be any Heegner point, and let $z$ be the divisor class:
$$
\frac{1}{\#(\OEx/\OFx)}\sum\limits_{\sigma\in\Gal(H/E)}\phi(x^\sigma),
$$
which is rational over $E$. This is the Heegner divisor in the
introduction, and its $p$-adic height will be related to the
derivative of the $p$-adic $L$ function.

\section{The geometric kernel}

\subsection{$p$-adic height pairing}

We are going to recall the basics of $p$-adic height pairings,
following Zarhin \cite{Zarhin90} and Nekov\'a\v{r}
\cite{Nekovar93}.\\

In this subsection, let $K$ be a number field and let $A/K$ be an
Abelian variety over $K$, which has good reduction at every place
$v$ of $K$ above $p$. For each integer $n\geq1$, write $S(A/K,n)$
for the $n$-Selmer group and write $S_p(A/K)$ for the limit
$\varprojlim\limits_{k}S(A/K,p^k)$. Also write
$T_p(A)=\varprojlim\limits_{k}A[p^k]$ for the $p$-adic Tate module,
and put $V_p(A)=T_p(A)\otimes_{\Zp}\Qp$.\\

The injection $A(K)\otimes\Zp\hookrightarrow S_p(A/K)$ gives, after
tensoring $\Qp$, an injection of $A(K)\otimes\Qp$ into
$S_p(A/K)\otimes_{\Zp}\Qp$, which, by \cite{BlochKato} 3.11, is
equal to $H_{f}^{1}(K,V_p(A))$, the Bloch-Kato Selmer group.\\

Let $K_{\infty}$ be the compositum of all $\Zp$-extension of $K$.
Suppose that for every place $v\mid p$, we are given a $\Qp$-linear
splitting for the Hodge Filtration
$\F^{0}H_{dR}^{1}(A/K_v)\rightarrow H_{dR}^{1}(A/K_v)$, which we fix
from now on. Then, by results of \cite{Zarhin90} and
\cite{Nekovar93}, a pairing
$$
H_{f}^1(K,V_p(A))\times
H_f^1(K,V_p(\hat{A}))\rightarrow\Gal(K_\infty/K)\otimes_{\Zp}\Qp
$$
is defined. This is the abstract $p$-adic pairing. If we are also
given a continuous morphism
$\ell:\A_{K}^{\times}/K^{\times}\rightarrow\Qp$, it will induce a
morphism from $\Gal(K_{\infty}/K)\otimes_{\Zp}\Qp$ to $\Qp$, so the
above pairing gives a pairing taking values in $\Qp$ after composing
with $\ell$.\\

\begin{rem}
We may replace $V_p(\hat{A})$ by $V_p(A)^{*}(1)$ via the Weil
pairing $V_p(A)\times V_p(\hat{A})\rightarrow\Q_p(1)$.
\end{rem}

Now we specialize to the case of Jacobians of curves and briefly
recall the definition of the height pairing. Let $X$ be an
irreducible smooth projective curve over $K$ with good reduction at
all places above $p$, and let $A$ be the Jacobian of $X$, which is
canonically isomorphic to its dual $\hat{A}$.\\

In this case, we have the $p$-adic Abel-Jacobi map: $A(K)\rightarrow
H^1(K,V_p(A))=\Ext^1_{\Gal(\overline{K}/K)}(\Qp,V_p(A))$. Thus for
each divisor $D$ of $X$ of degree $0$, we may represent the image of
$D$ under the Abel-Jacobi map by an extension:
$$
0\rightarrow
V_p(A)\rightarrow\mathcal{E}_D\rightarrow\Qp\rightarrow0.
$$

Let $D_1$, $D_2$ be two divisors of $X$ of degree $0$ which do not
intersect. The height pairing $\langle cl(D_1),cl(D_2)\rangle$ is a
sum of local pairings (here $cl$ means the class in $A(K)$). The
decomposition depends on a mixed extension $\mathcal{E}$, i.e. a
$\Gal(\overline{K}/K)$-module $\mathcal{E}$ that fits into the
commutative diagram:
$$
\xymatrix{
& & 0\ar[d] & 0 \ar[d]\\
0\ar[r] & \Qp(1)\ar[r]\ar[d]^= &
\mathcal{E}_{D_2}^{*}(1)\ar[r]\ar[d] &
V_p(A)\ar[r]\ar[d] & 0\\
0\ar[r] & \Qp(1)\ar[r] & \mathcal{E}\ar[r]\ar[d] &
\mathcal{E}_{D_1}\ar[r]\ar[d] & 0\\
& & \Qp\ar[r]^=\ar[d] & \Qp\ar[d]\\
& & 0 & 0}
$$
Such a mixed extension arise naturally as subquotient of the
relative cohomology
$H_{et}^1(X_{\overline{K}}-\left|D_1\right|,\left|D_2\right|;\Qp)(1)$,
c.f. \cite{Nekovar95} II 1.9.\\

If $v$ is a place not dividing $p$, then we have
$H^1(K_v,V_p(A))=H^1(K_v,V_p(A)^{*}(1))=0$, thus as a
$\Gal(\overline{K_v}/K_v)$-module, the module $\mathcal{E}$ splits
into a direct sum $V_p(A)\oplus\mathcal{E}_v'$, where
$\mathcal{E}_v'$ is an extension of $\Qp(1)$ by $\Qp$. We denote by
$\mathcal{C}_v$ the class of $\mathcal{E}_v'$ in the group
$\Ext^1_{\Gal(\overline{K_v}/K_v)}(\Qp,\Qp(1))$.\\

On the other hand, we have:
$$
\Ext^1_{\Gal(\overline{K_v}/K_v)}(\Qp,\Qp(1))=
H^1(K_v,\Qp(1))=\left(\varprojlim\limits_{k}K_v^{\times}/
(K_v^{\times})^{p^k}\right)\otimes_{\Zp}\Qp=K_v^{\times}
\widehat{\otimes}\Qp.
$$
Let $\ell_v$ be the map
$K_v^{\times}\rightarrow\A_K^{\times}/K^{\times}
\underrightarrow{\ell}\Qp$, where $\ell$ is a fixed continuous
morphism from $\A_K^{\times}/K^{\times}$ to $\Qp$, then the local
height pairing is defined as $\langle D_1,
D_2\rangle_v=-(\ell_v\otimes1)(\mathcal{C}_v)$.\\

If $v$ is a place dividing $p$, one constructs similarly a class
$\mathcal{C}_v$ in the group $H^1(K_v,\Qp(1))$, and defines the
local pairing by the same formula. For details, see \cite{Nekovar95}
II 1.7 and \cite{Nekovar93}.\\

In the following, we always take the morphism $\ell$ to be the
logarithm of the cyclotomic character: $\ell=\log_p\circ\xi_p$.

\subsection{Construction of the geometric kernel}
We will define the geometric kernel $\Psi$ to be a $p$-adic modular
form in the space $\MS_{2}(N,L)$ with $q$-expansion given by:
$a(I,\Psi)=\langle z,\TN(I)z\rangle$, where $z$ is the Heegner
divisor in subsection \ref{heegnerdivisor}. By duality of $p$-adic
modular forms (c.f. \cite{Hida91} section 2 for a setting close to
ours), this $q$-expansion is really a $p$-adic
modular form.\\

We also define:
\begin{eqnarray*}
\Psi_{p}=\sum\limits_{I}\sum\limits_{w\mid p}\langle
z,\TN(I)z\rangle_{w}q^{I}\\
\Psi_{f}=\sum\limits_{I}\sum\limits_{w\nmid p}\langle
z,\TN(I)z\rangle_{w}q^{I},
\end{eqnarray*}
which are also in the space $\MS_{2}(N,L)$ and satisfy
$\Psi_{p}+\Psi_{f}=\Psi$.

\section{$p$-adic heights at places outside $p$}

In this section, let $v$ be a finite place of $F$ not dividing $p$.\\

In this case the $p$-adic height pairing
$\langle\cdot,\cdot\rangle_{w}$ is unique, and takes values in the
subset $\Q\cdot\log_{p}(\left|q_{v}\right|_\infty)$. Since the
corresponding real height pairings
$\langle\cdot,\cdot\rangle_{w,\infty}$ takes values in the subset
$\Q\cdot\log(\left|q_{v}\right|_\infty))$ of $\R$ and satisfies the
same conditions as the $p$-adic one, by uniqueness we have:
$$
\langle\cdot,\cdot\rangle_{w}/\log_{p}(\left|q_{v}\right|_\infty))=
\langle\cdot,\cdot\rangle_{w,\infty}/\log(\left|q_{v}\right|_\infty)).
$$
In \cite{Zhang01}, the real height pairings are computed, and we are
going to use these results to show that the form $\Psi_{f}$ is
closely related to the analytic kernel $\Phi$.

\subsection{The quotient space}

Following \cite{Zhang01}, we define $S$ to be the space of functions
on the set of non-zero integral ideals of $\OF$ with values in the
$p$-adic field $L$, modulo the equivalence relation that two
functions $a$ and $b$ are equivalent if there is a non-zero ideal
$M$ such that $a(I)=b(I)$ for all $I$ prime to $M$. Identifying a
$p$-adic modular form with its $q$-expansion, we may view the space
$\MS_{2}(N,L)$ as a subspace of $S$.\\

A function $f$ in $S$ is called "quasi-multiplicative" if there is a
non-zero ideal $M$ such that $f(IJ)=f(I)f(J)$ for all $I$, $J$ such
that $(I,J)=(M,IJ)=1$. If $f$ is quasi-multiplicative function, a
function $h$ is called an "$f$-derivative" if
$h(IJ)=f(I)h(J)+h(I)f(J)$ for all $I$, $J$ as above.\\

The functions $I\mapsto\sigma_1(I)$ and $I\mapsto r(I)$ define
functions $\sigma_1$ and $r$ in the space $S$. Let $D_N$ be the
subspace of $S$ generated by $\sigma_1$, $r$,
$\sigma_1$-derivatives, $r$-derivatives, and all forms that are "old
at $N$", i.e. comes from a form of level $Mp^{\infty}$ with $M\mid
N$ and $M\neq N$.\\

\begin{prop}
If the $q$-expansion of a form in $\MS_{2}(N,L)$ lies in $D_N$, then
it is old at $N$.
\end{prop}

\begin{dem}
This is the same as Proposition 4.5.1 of \cite{Zhang01}.
\end{dem}

In view of this proposition, forms that differ by some function in
$D_N$ will have the same value under the linear functional $l_{f}$.

\subsection{Compare the two kernels}

In \cite{Zhang01}, the real heights of the Heegner divisors are
calculated, and we may deduce from it the corresponding $p$-adic
heights, i.e. coefficients of the $q$-expansion of the geometric
kernel $\Psi_f$.

\begin{prop}
Modulo $D_N$, the function $I\mapsto a(I,\Psi_f)$ on the set of
non-zero ideals of $\OF$ is equal to a sum
$\sum\limits_{v}a_v(I,\Psi_f)$, where the sum ranges over all finite
places of $F$, and each $a_{v}(I,\Psi_f)$ is given by:
\begin{enumerate}
\item if $v$ is inert in $E$ and $q_v$ is the associated prime
ideal, then:
$$
a_v(I,\Psi_f)=\sum\limits_{\substack{\alpha+\beta=1,\alpha,\beta>0\\
\epsilon(-\alpha\beta)=1,\forall w\mid D\\N\mid\beta DI}}\delta
r(\alpha I)r(\beta I/Nq_v)\ord_v(\beta
Iq_v/N)\log_p(\left|q_{v}\right|_\infty),
$$
here $\delta$ and $r$ has the same meaning as in Proposition
\ref{analyticcoeff};
\item if $v$ is ramified in $E$ and $q_v$ is the associated prime
ideal, then:
$$
a_v(I,\Psi_f)=\sum\limits_{\substack{\alpha+\beta=1,\alpha,\beta>0\\
\epsilon_{w}(-\alpha\beta)=1,\forall v\neq w\mid D\\N\mid\beta
DI\\\epsilon_v(-\alpha\beta)=-1}}\delta r(\alpha I)r(\beta
I/N)\ord_v(\beta Iq_v)\log_p(\left|q_{v}\right|_\infty);
$$
\item if $v$ is split in $E$, then $a_v(I,\Psi_f)=0$.
\end{enumerate}
\end{prop}

\begin{dem}
This is the combination of \cite{Zhang01} Proposition 5.4.8,
Proposition 6.4.5 and Proposition 7.1.1, but we changed the real
logarithm to the $p$-adic logarithm everywhere.
\end{dem}

The relation between the two forms $\Phi$ and $\Psi_f$ can then be
stated as follows.

\begin{prop}
The difference of the two forms $\left(\prod\limits_{\PP_i\mid
p}(\T(\PP_i)^4-\T(\PP_i)^2)\right)\Phi$ and
$\left(\prod\limits_{\PP_i\mid p}(\T(\PP_i)-1)^4\right)\Psi_f$ is
annihilated by the linear functional $l_f$.
\end{prop}

\begin{dem}
This is essentially the same as \cite{Perrin-Riou87} Proposition
3.20, but note that our $\Psi$ takes the sum of the heights of all
conjugates of $x$, thus the translation by an element
$\sigma\in\Gal(H/E)$ acts trivially.
\end{dem}

If we can prove that the form $\left(\prod\limits_{\PP_i\mid
p}(\T(\PP_i)-1)^4\right)\Psi_p$ is also annihilated by the linear
form $l_f$, then applying $l_f$ will give:
\begin{eqnarray*}
& & \left(\prod\limits_{i}(\alpha_i^4-\alpha_i^2)\right)
L_{p,\xi_E}'(1)=l_f\left(\left(\prod\limits_{i}
(\T(\PP_i)^4-\T(\PP_i)^2)\right)\Phi\right)\\
&=& l_f\left(\left(\prod\limits_{i}(\T(\PP_i)-1)^4\right)\Psi\right)
=\left(\prod\limits_{i}(\alpha_i-1)^4(1-\frac{\left|\PP_i
\right|_\infty}{\alpha_i^2})^{-1}\right)\frac{(f,\Psi)}{(f,f)},
\end{eqnarray*}
where the last equality comes from Lemma \ref{lemmaproduct}. Since
we have $\frac{(f,\Psi)}{(f,f)}=\langle z_f,z\rangle=\langle
z_f,z_f\rangle$, this gives our main result.\\

Thus it only remains to prove that the contribution of the form
$\Psi_p$ is zero. Actually, we will prove that for each place $v$ of
$F$ over $p$, the contribution of the form
$\Psi_v=\sum\limits_{I}\sum\limits_{w\mid v}\langle
z,\T(I)z\rangle_{w}q^I$ is zero.

\section{$p$-adic heights at places above $p$}

In this section we prove that the contribution of the local heights
for places above $p$ is zero. We are going to use a corrected
version of the method in \cite{Nekovar95}.\\

Thus recall that $X/F$ is the Shimura curve constructed in
subsection \ref{curveX} and $A$ is the Jacobian of $X$.

\subsection{General arguments}

We first decompose the space of modular forms $\MS_2(N)$ according
to action of the spherical Hecke algebra. Let $K$ be a number field,
big enough to contain all Hecke eigenvalues of eigenforms in
$\MS_2(N)$. The spherical Hecke algebra
$\mathbb{T}=K[\TN(I),(I,Np)=1]$, a subalgebra of
$\End_{L_0}(\MS_2(N,K))$, then decomposes as a direct sum of
$\mathbb{T}$ modules:
$$
\mathbb{T}=\bigoplus\limits_{\rho:\mathbb{T}\rightarrow K}K,
$$
where $\rho$ runs through all morphisms from $\mathbb{T}$ to $K$.
For each $\rho$ as above, let $e_\rho$ be the projector onto the
$\rho$-component.

\begin{lem}\label{spherical}
There is a finite set $S$ of ideals $I$ of $\OF$ satisfying
$(I,Np)=1$ and $r(I)=0$, such that every projector $e_\rho$ is a
linear combination of the Hecke operators $\TN(I)$.
\end{lem}

\begin{dem}
This follows from \cite{Nekovar95} Lemma II.5.7.
\end{dem}

Let $K_p$ be the closure of $K$ in $\Cp$, under the fixed embedding
$\overline{\Q}\hookrightarrow\overline{\Qp}$. We then have the
corresponding decomposition of $V_p(A)$ as representation of
$\Gal(\overline{\Q}/F)$ over $K_p$:
$$
V_p(A)\otimes_{\Qp}K_p=\bigoplus\limits_{\rho:\mathbb{T}\rightarrow
K}V_\rho(1)^{m_\rho},
$$
where $V_\rho$ is a $K_p$ representation of $\Gal(\overline{\Q}/F)$
of dimension $2$, characterized by the identity:
$$
\det(1-X\cdot\Frob(Q)\big|_{V_\rho})=1-\rho(\TN(Q))X+\left|Q\right|
_{\infty}X^2
$$
for any prime ideal $Q$ prime to $Np$, and $m_\rho\geq1$ is the
multiplicity.\\

Let $\rho_f:\mathbb{T}\rightarrow K$ be the morphism corresponding
to our fixed form $f$. Since $f$ is a new form, the multiplicity
$m_{\rho_f}$ is equal to $1$, by \cite{Zhang01}, Theorem 3.2.1.
Hence the space $V=V_p(A)\otimes_{\Qp}K_p$ decomposes as:
$V=V'\oplus V''$, with $V'=e_{\rho_f}V$ a representation of
dimension $2$, the $f$-component, and
$V''=(1-e_{\rho_f})V$.\\

Now fix a prime ideal $\PP_i$ of $\OF$ above $p$. If we view $V'$ as
representation of the Galois group $\Gal(\overline{\Qp}/F_{\PP_i})$,
then the "reduction mod $\PP_i$" map on the variety $A$ gives a
short exact sequence:
$$
0\rightarrow V'_{+,i}\rightarrow V'\rightarrow V'_{-,i}\rightarrow0,
$$
with $V'_{-,i}\simeq(V'_{+,i})^{*}(1)$ an unramified representation,
on which the Frobenius acts as multiplication by a Weil number of
weight $-1$.\\

Let $T'$ and $T''$ be fixed $\Gal(\overline{\Q}/\Q)$-stable
$\OF_{K_p}$-lattices of $V'$ and $V''$, respectively. Put
$T'_{+,i}=T'\cap V'_{+,i}$ and $T'_{-,i}=T'/T'_{+,i}$.

\begin{prop}\label{prop1}
If $\Omega$ is an algebraic extension of $F_{\PP_i}$, such that the
maximal unramified subextension $\Omega^{ur}$ is of finite degree
over $F_{\PP_i}$, then we have $H^0(\Omega,V'_{-,i})=0$.
\end{prop}

\begin{dem}
The representation $V'_{-,i}$ being unramified, the non-triviality
of $H^0(\Omega,V'_{-,i})$ would imply non-triviality of the
cohomology $H^0(\Omega^{ur},V'_{-,i})$, hence implies that some
power of the Frobenius morphism acts trivially on the space
$V'_{-,i}$. But none of the powers of a Weil number is equal to $1$,
a contradiction.
\end{dem}

\begin{cor}\label{cor1}
The groups $H^0(\Omega,V'_{-,i}/T'_{-,i})$ and
$H^0(\Omega,(V'_{+,i})^{*}(1)/(T'_{+,i})^{*}(1))$ are finite.
\end{cor}

\begin{dem}
Temporarily write $V$ for $V'_{-,i}$, $T$ for $T'_{-,i}$ and $G$ for
$\Gal(\overline{\Qp},\Omega)$. If the group
$H^0(\Omega,V'_{-,i}/T'_{-,i})$ is infinite, then there is a
sequence $(v_m)_m$ in $V$, all different modulo $T$, such that
$v_m^{\sigma}-v_m$ is in $T$ for every $\sigma\in G$. Since the set
$p^{-k}T/T$ is finite for every $k$, the sequence $(v_m)_m$ is
unbounded, i.e. for every $k$, there is a term $v_{m_k}$ not
belonging to $p^{-k}T$. Fix such $v_{m_k}$ for every $k$.\\

Let $a_k(\geq k)$ be the integer such that $v_{m_k}\in
p^{-a_k-1}T\backslash p^{-a_k}T$, and put $u_k=p^{a_k}v_{m_k}$,
which lies in $p^{-1}T\backslash T$. We then have
$u_k^{\sigma}-u_k\in p^{a_k}T\subseteq p^{k}T$ for any $\sigma\in
G$. Since the set $p^{-1}T\backslash T$ is compact, the sequence
$(u_k)_k$ admits at least one accumulation point $u$, which must be
a non-zero fixed point of $G$, contradicting the above
proposition.\\

Changing the meanings of $V$ and $T$ proves that
$H^0(\Omega,(V'_{+,i})^{*}(1)/(T'_{+,i})^{*}(1))$ is finite.
\end{dem}

Now let $T$ be the lattice $T_p(A)\otimes_{\Zp}\OF_{K_p}$ of
$V=V_p(A)\otimes_{\Qp}K_p$. Take $T'=e_{\rho_f}T$ and
$T''=(1-e_{\rho_f})T$, which are lattices of $V'$ and $V''$,
respectively. Then the sum $T'\oplus T''$ is again a lattice of $V$,
thus there is a constant $d_0$, such that: $T\subseteq T'\oplus
T''\subseteq p^{-d_0}T$.\\

We will also fix a splitting for the Hodge filtration
$F^0(H_{dR}^1(A/F_{\PP_i})\otimes_{\Qp}K_p)\hookrightarrow
H_{dR}^1(A/F_{\PP_i})\otimes_{\Qp}K_p$. Under the decomposition
$$
H_{dR}^1(A/F_{\PP_i})\otimes_{\Qp}K_p=D_{cris}(V'\big|
_{\Gal(\overline{\Qp}/F_{\PP_i})})\{-1\}\oplus
D_{cris}(V''\big|_{\Gal(\overline{\Qp}/F_{\PP_i})})\{-1\},
$$
the $F^0$ part decomposes as: $D_{cris}(V'_{-,i})\{-1\}\oplus
(F^{1}D_{cris}(V''\big|_{\Gal(\overline{\Qp}/F_{\PP_i})}))\{-1\}$.
We fix the canonical splitting on the first factor, and any
splitting on the second factor.\\

\begin{rem}
If the curve have good ordinary reduction at $\PP_i$ (which is our
case), the $p$-adic height pairing given by this splitting is the
same as that defined by the universal norm, as in \cite{Nekovar93}
Theorem 6.11.
\end{rem}

We fix a proper smooth model $\mathcal{X}_{i}$ of $X$ over the ring
$\OF_{\PP_i}$. Such a model exists under our hypothesis.

\begin{prop}\label{prop2}
Let $\Omega$ be a fixed extension of $F_{\PP_i}$ as in Proposition
\ref{prop1}, and let $H$ be a finite subextension of $\Omega$. If
$D_1$ and $D_2$ are two divisors of $X$ of degree zero over $H$,
such that the Zariski closures of $D_1$ and $e_{\rho_f}(D_2)$ in
$\mathcal{X}_i$ do not intersect, then there is a constant $C$,
depending only on $\Omega$, such that:
$$
\langle D_1,p^{d_0}e_{\rho_f}(D_2)\rangle_H\in
p^{-C}\log_p(\xi_H(\OF_{H}^{*}\widehat{\otimes}\Zp)).
$$
\end{prop}

\begin{dem}
The proof is in \cite{Nekovar95} II Proposition 1.11, with the
following modifications:\\

1. By definition of $d_0$, we see that $p^{d_0}e_{\rho_f}(D_2)$ is
in $T'$, so there is no $T''$ component, thus the constant $d_1$ in
\textit{loc. cit.} can be replaced by $0$;\\

2. In this case the mixed extension $\mathcal{E}$ is crystalline as
a $\Gal(\overline{\Qp}/H)$-representation, thanks to the comparison
theorem for relative cohomology (\cite{Olsson05} Theorem 13.21), so
we can replace $H^{*}\widehat{\otimes}\Zp$ by
$\OF_{H}^{*}\widehat{\otimes}\Zp$;\\

3. By \cite{Nekovar93} 6.9, the constant $p^{d_2}$ in \textit{loc.
cit.} divides the number
$$
\#H^0(H,V'_{-,i}/T'_{-,i})\#H^0(H,(V'_{+,i})^*(1)/(T'_{+,i})^*(1)),
$$
which is bounded by the same number with the field $H$ replaced by
$\Omega$. This last number is finite by Corollary \ref{cor1}, and
only depends on the field $\Omega$.
\end{dem}

\subsection{Application to Heegner divisors}

We now apply the above general arguments to Heegner divisors on the
curve $X$ with complex points
$X(\C)=\Bx\backslash\mathcal{H}^{\pm}\times\widehat{B}^{\times}/
\widehat{F}^{\times}\widehat{R}^{\times}$. If $b$ is any element of
$\widehat{B}^{\times}$, write $c(b)$ for the CM point represented by
$(z_0,b)$, with $z_0$ the unique fixed point of $\Ex$ on the upper
half plane.\\

Fix throughout this subsection a prime ideal $\PP_i$ of $\OF$ above
$p$, and a place $v$ of the conductor $1$ ring class field $H$ over
$\PP_i$. Let $\PP_i=q_i\overline{q}_i$ be the decomposition of
$\PP_i$ in $\OEE$.\\

Let $x$ be a Heegner point as in subsection \ref{heegnerdivisor}. If
$I$ is an ideal of $\OF$ such that $r(I)=0$, then the divisors $x$
and $\TN(I)(x)$ do not intersect. By Serre-Tate theory (c.f.
\cite{Zhang01} section 5.2.2, Case 1), their Zariski closures on the
integral model $\mathcal{X}_i$ are also
disjoint.\\

Let $\xi$ be a divisor representing the Hodge class which does not
intersect with any $\TN(I)(x)$. For each prime divisor $D$ occurring
in $\xi$, let $Z_D$ be the set of ideals $I$ of $\OF$ such that
$(I,N)=1$ and the Zariski closure $\mathcal{D}$ of $D$ on
$\mathcal{X}_i$ intersects the Zariski closure of $\TN(I)(x)$.\\

\begin{lem}
If $I$ is an ideal in $Z_D$ and $J$ is another ideal such that
$r(J)=0$, then the ideal $IJ$ is not in $Z_D$.
\end{lem}

\begin{dem}
The hypothesis $r(J)=0$ implies that the divisors $\TN(I)(x)$ and
$\TN(IJ)(x)$ do not intersect, thus their Zariski closures on
$\mathcal{X}_i$ do not intersect, again by \cite{Zhang01} section
5.2.2, Case 1. Since the intersection of $\mathcal{D}$ and the
special fiber is a prime divisor, this implies that the Zariski
closure of $\TN(IJ)(x)$ does not intersect with $\mathcal{D}$ on the
special fiber.
\end{dem}

Let $M$ be the ideal $\prod\limits_{Q\in S}Q$, where $S$ is the
finite set in Lemma \ref{spherical}. For any prime divisor $D$ in
the support of $\xi$ such that the set $Z_D$ is non-empty, choose an
ideal $I_D$ of $Z_D$. Let $\mathcal{I}$ be the product
$\prod\limits_{D}I_D$. Let $J$ be any ideal of $\OF$ prime to
$\mathcal{I}NM$, and let $\mathcal{J}$ be the ideal $\mathcal{I}J$.
The above lemma then shows that the Zariski closures of $\xi$ and
$\TN(I)(x)$ do not intersect if $I$ is an ideal prime to $N$ and
divisible by $\mathcal{J}$.\\

\begin{rem}
The ideal $\mathcal{J}$ depends on $i$, but in the following we may
replace it with the product for all $i$, thus the above property
will be valid for every model $\mathcal{X}_i$.
\end{rem}

We choose another representative $\xi'$ of the Hodge class, such
that its Zariski closure on $\mathcal{X}_i$ does not intersect that
of $x-\xi$. The moving lemma on the special fiber of $\mathcal{X}_i$
guarantees that such a divisor exists for every $i$, and then
applying Chinese remainder theorem gives a divisor $\xi'$ that
satisfies the above property for each $i$.\\

For every element $\sigma\in\Gal(H/E)$, put:
\begin{eqnarray*}
\Psi_{v,\sigma} &=& \sum\limits_{I}\langle
x-\xi,\TN(I)(x^{\sigma}-\xi')\rangle_{v}q^I\\
&=& \sum\limits_{I}\langle
x-\xi,\TN(I)(x^{\sigma})-\deg(\TN(I))\xi'\rangle_{v}q^I
\end{eqnarray*}
as a formal $q$-expansion. We let an element $\tau\in\Gal(H/E)$ act
on functions $A$ over the group $\Gal(H/E)$ by translation:
$$
\tau A(\sigma)=A(\tau\sigma).
$$
Write $\sigma_{q_i}$ for the Frobenius element of $\Gal(H/E)$
corresponding to the ideal $q_i$, and put:
\begin{eqnarray*}
U_{v,\sigma} &=& (\T(\PP_i)-\sigma_{q_i})(\sigma_{q_i}\T(\PP_i)-1)
\Psi_{v,\sigma}\\
&=& \sum\limits_{I}\langle
x,(\TN(I\PP_i^2)+\TN(I))(x^{\sigma\sigma_{q_i}})-\TN(I\PP_i)
(x^{\sigma}+x^{\sigma\sigma_{q_i}^2})\rangle_v.
\end{eqnarray*}

\begin{lem}\label{normrelation}
("norm relation", or "trace relation") Let $k\geq0$ be an integer.
There is a divisor $D_k$ defined over the ring class field
$H[\PP_i^{k+2}]$, such that the trace $\tr_{H[\PP_i^{k+2}]/H}(D_k)$
is equal to the divisor
$(\TN(I_k\PP_i^{2})+\TN(I_k))(x^{\sigma\sigma_{q_i}})-\TN(I_k\PP_i)
(x^{\sigma}+x^{\sigma\sigma_{q_i}^2})$, with $I_k:=I\PP_i^k$.
\end{lem}

We will prove this lemma in subsection \ref{sectionNormRel}. Here we
give the consequence. Let $d_0$ be the fixed constant as in
Proposition \ref{prop2}.

\begin{prop}
Let $I$ be an ideal of $\OF$ prime to $NM$ and divisible by
$\mathcal{J}$ such that $r(I)=0$. There is a constant $C'$, not
depending on $k$, such that the height $\langle
x,p^{d_0}e_{\rho_f}(\TN(I\PP_i^2)+\TN(I))(x^{\sigma\sigma_{q_i}})
-\TN(I\PP_i)(x^{\sigma}+x^{\sigma\sigma_{q_i}^2})\rangle_v$ has
value in $p^{k-C'}\Zp$.
\end{prop}

\begin{dem}
Let $w=w_k$ be a place of $H[\PP_i^{k+2}]$ above $v$. By Lemma
\ref{normrelation}, the above height is equal to: $\langle
x-\xi,p^{d_0}e_{\rho_f}D_k\rangle_{w}$. Since the divisors do not
intersect on the integral model, we may apply Proposition
\ref{prop2} with $\Omega$ equals to the union of all the extensions
$H[\PP_i^k]_{w_k}$ for all $k$. This gives:
$$
\langle x-\xi,p^{d_0}e_{\rho_f}D_k\rangle_{w}\in
p^{-C}\log_p(\xi_p(\N_{H[\PP_i^{k+2}]_w/H_v}
(\OF_{H[\PP_i^{k+2}]_w}^{*})\widehat{\otimes}\Zp)).
$$
The extension $H[\PP_i^{k+2}]_w/H_v$ is ramified with ramification
index $\N(\PP_i)^{k+2}$, so by class field theory we know that the
number
$\log_p(\xi_p(\N_{H[\PP_i^{k+2}]_w/H_v}(\OF_{H[\PP_i^{k+2}]_w}^{*})
\widehat{\otimes}\Zp))$ lives in $p^{k+2}\Zp$.
\end{dem}

Now take the same of all $\sigma\in\Gal(H/E)$ and all $x$ of
conductor $1$ in the above proposition, we get the similar result
with $x$ and $x^{\sigma}$ replaced by the divisor $z$ in subsection
\ref{heegnerdivisor}. Note that, since the translation by any
$\sigma\in\Gal(H/E)$ permutes the group, the operator
$(\T(\PP_i)-\sigma_{q_i})(\sigma_{q_i}\T(\PP_i)-1)$ is the same as
$(\T(\PP_i)-1)^2$ when acting on the sum. Also note that the sum of
$\Psi_{v}=\sum\limits_I\langle z,\T_N(I)z\rangle_vq^I$ for all
places $v$ is equal to the form $\Psi_p$.

\begin{cor}
We have $l_f\left(\prod\limits_{i}(\T(\PP_i)-1)^4\right)\Psi_p=0$.
\end{cor}

\begin{dem}
We have shown that there is a $p$-adic modular form $\Psi'$ in the
space $\MS_2(N,L)$ which has the same coefficients as
$\left(\prod\limits_{i}(\T(\PP_i)-1)^4\right)\Psi_p$ on ideals $I$
divisible by $\PP\mathcal{J}$ and such that $(I,NM)=1$ and $r(I)=0$.
For these $I$, we have $a(I,e_{\rho_f}e_\ord(\Psi'))=0$ by letting
$k$ tend to infinity in the above proposition. This implies that the
form $e_{\rho_f}e_\ord(\Psi')$ is equal to zero, by a variant of
\cite{Nekovar95} Lemma 5.7.
\end{dem}

\begin{rem}
The remark in \cite{Howard05} Remark 2.0.3 is not valid, because the
form $e_{\rho_f}e_\ord(\Psi')$ is a multiple of the form $f_0$,
hence the counterexample in \textit{loc. cit.} does not apply.
\end{rem}

\subsection{The norm relation}\label{sectionNormRel}

In this subsection we prove Lemma \ref{normrelation}. Without loss
of generality, we may assume that the ideal $I$ has trivial
prime-to-$\PP_i$ part, and that the Heegner point $x^{\sigma}$ is
the "basic Heegner point" represented by $(z_0,1)$, with the complex
description as in subsection \ref{sectionCM}.\\

Let $CM(X)$ denote the set of CM points on $X$. Since the prime
ideal $\PP_i$ splits in $\OEE$ as $\PP_i=q_i\overline{q}_i$, the
algebra $E_{\PP_i}=E\otimes_{F}F_{\PP_i}=E_{q_i}\times
E_{\overline{q}_i}$ is canonically isomorphic to $F_{\PP_i}^2$. We
fix an isomorphism $B\otimes_{F}F_{\PP_i}\simeq M_2(F_{\PP_i})$,
such that the restriction to the order $R\otimes_{\OF}\OF_{\PP_i}$
gives an isomorphism to $M_2(\OF_{\PP_i})$, and that the image of
the element $(a,b)\in F_{\PP_i}^2\simeq E_{\PP_i}$ is the diagonal
matrix $
\begin{pmatrix}
a & \\
 & b
\end{pmatrix}$.\\

Once such an isomorphism is fixed, we get a map:
\begin{eqnarray*}
c:\GL(F_{\PP_i})\simeq(B\otimes_{F}F_{\PP_i})^{\times} &\rightarrow&
CM(X)\\
b &\mapsto& [z_0,b],
\end{eqnarray*}
here $[z_0,b]$ stands for the CM point represented by the pair
$(z_0,b)$. This map obviously factors through the quotient:
$\PGL(F_{\PP_i})/\PGL(\OF_{\PP_i})=\mathcal{T}$, the set of vertices
of the Bruhat-Tits tree. We write $d(\cdot,\cdot)$ for the distance
function on this tree.\\

There is a special line $\mathcal{L}$ on the tree $\mathcal{T}$,
given by the image of $E_{\PP_i}^{\times}$, or equivalently, the
image of matrices of the form $
\begin{pmatrix}
F_{\PP_i}^{*} & \\
 & F_{\PP_i}^{*}
\end{pmatrix}$. If $b$ is a vertex on the line $\mathcal{L}$, then
by Shimura's reciprocity law, the corresponding CM point $c(b)$ is
defined over the conductor $1$ ring class field $H$, because $b$
commutes with any element $\sigma\in\widehat{E}^{\times}$.
Similarly, if $b$ is a vertex of distance $k$ to the line
$\mathcal{L}$, then the stabilizer of $b$ in $\widehat{E}^{\times}$
contains the $\PP_i^k$-th order, hence the point $c(b)$ is define
over the ring class field $H[\PP_i^k]$. The Galois group
$\Gal(H[\PP_i^k]/H)$ acts as isometry on the tree and fixes the line
$\mathcal{L}$, and for a given point $x$ on the line, the action
permutes all the
points of distance $k$ to the point $x$ that are not on the line.\\

Let $\pi_i\in F_{\PP_i}$ be a uniformiser. The Hecke operators (or
correspondences) $\TN(\PP_i^k)$ act in the usual way on the tree
$\mathcal{T}$:
$$
\TN(\PP_i^k)(x)=\sum\limits_{\substack{d(x,y)\leq k\\d(x,y)\equiv
k\mod2}}y.
$$
It is then easy to see that, under the map $c$, these actions
coincide with the actions of the Hecke operators on the divisors of
$X$. (More precisely, this follows from the following decomposition:
$$
U(\PP_i^k)=K
\begin{pmatrix}
\pi_i^k & \\
 & 1
\end{pmatrix}K\bigsqcup K
\begin{pmatrix}
\pi_i^{k-1} & \\
& \pi_i
\end{pmatrix}K\bigsqcup\cdots\bigsqcup K
\begin{pmatrix}
\pi_i^{\lceil k/2\rceil} & \\
 & \pi_i^{\lfloor k/2\rfloor}
\end{pmatrix}K,
$$
where $U(\PP_i^k)$ is the set of matrices $g$ in $M_2(\OF_{\PP_i})$
such that $(\det g)\OF_{\PP_i}=\PP_i^k$, which was used to define
the Hecke operators, and $K$ is the group $\GL(\OF_{\PP_i})$.)\\

We also have the action of the Frobenius morphism $\sigma_{q_i}$,
which sends the basic Heegner point $x$ to the point $c
\begin{pmatrix}
\pi_i & \\
 & 1
\end{pmatrix}$.\\

Now we can prove the lemma using the description of the tree. Put
$x_A=x^{\sigma}$, $x_B=x^{\sigma\sigma_{q_i}}$,
$x_C=x^{\sigma\sigma_{q_i}^2}$, these are three consecutive points
on the line $\mathcal{L}$. We then have:
\begin{eqnarray*}
&& (\TN(\PP_i^{k+2})+\TN(\PP_i^k))(x_B)-\TN(\PP_i^{k+1})(x_A+x_C)\\
&=& \sum\limits_{\substack{d(y,x_B)\leq k+2\\d(y,x_B)\equiv
k+2\mod2}}y-\sum\limits_{\substack{d(y,x_A)\leq k+1\\d(y,x_A)\equiv
k+1\mod2}}y-\sum\limits_{\substack{d(y,x_C)\leq k+1\\d(y,x_C)\equiv
k+1\mod2}}y+\sum\limits_{\substack{d(y,x_B)\leq k\\d(y,x_B)\equiv
k\mod2}}y.
\end{eqnarray*}
Write $R_1$, $R_2$, $R_3$ and $R_4$ for the four sets of summation
in the above formula. It is then clear that $R_2$, $R_3$ and $R_4$
are all subsets of $R_1$, and that the intersection of $R_2$ and
$R_3$ is exactly $R_4$. The inclusion-exclusion principle then tells
us that the above is equal to the sum over the set
$R_1\backslash(R_2\bigcup R_3)$, which is the set of vertices $y$
such that $d(y,x_B)=k+2$ and that the shortest path between $y$ and
$x_B$ does not pass any other vertex on the line $\mathcal{L}$. Let
$y$ be any one of these vertices, applying the map $c$ to the above
formula shows that the divisor in Lemma \ref{normrelation} is equal
to the trace $\tr_{H[\PP_i^{k+2}]/H}(c(y))$.

\end{document}